\documentclass{amsart}
\usepackage{amsmath}
\usepackage{amscd}
\usepackage{amsfonts,latexsym,amssymb}

\renewcommand{\a}{\alpha}
\renewcommand{\b}{\beta}
\renewcommand{\d}{\delta}
\newcommand{\g}{\gamma}
\newcommand{\e}{\varepsilon}
\newcommand{\f}{\epsilon}

\renewcommand{\l}{\lambda}
\newcommand{\s}{\sigma}

\renewcommand{\o}{\omega}
\newcommand{\p}{\phi}

\newcommand{\D}{\Delta}

\renewcommand{\P}{\Phi}

\newcommand{\SE}{{\mathcal{E}}}

\newcommand{\SU}{{\mathcal{U}}}

\DeclareFontFamily{OT1}{rsfs}{}
\DeclareFontShape{OT1}{rsfs}{n}{it}{<->rsfs10}{}
\DeclareMathAlphabet{\curly}{OT1}{rsfs}{n}{it}

\newcommand{\CrP}{{\curly P}}

\newcommand{\CrV}{{\curly V}}

\newcommand{\PP}{\mathbb{P}}

\newcommand{\C}{\mathbb{C}}

\newcommand{\N}{\mathbb{N}}
\newcommand{\R}{\mathbb{R}}

\newcommand{\CP}{\mathbb{CP}}

\newcommand{\Hom}{\operatorname{Hom}}

\newcommand{\id}{\textbf{\textit{I}}}
\newcommand{\im}{\operatorname{im}}
\newcommand{\dimC}{\text{dim}_{\C}\,}
\newcommand{\codimC}{\text{codim}_{\C}\,}
\newcommand{\surj}{\twoheadrightarrow}
\newcommand{\inj}{\hookrightarrow}

\newcommand{\too}{\longrightarrow}
\newcommand{\pr}{\operatorname{pr}}
\newcommand{\End}{\operatorname{End}}
\newcommand{\Tr}{\operatorname{Tr}}
\newcommand{\PD}{\operatorname{PD}}

\newcommand{\la}{\langle}
\newcommand{\ra}{\rangle}

\newcommand{\sref}{s_{k,x}^{\text{\scriptsize ref}}}
\newcommand{\bd}{\partial}
\newcommand{\bbd}{\bar{\partial}}
\newcommand{\x}{\times}
\newcommand{\ox}{\otimes}

\newtheorem{proposition}{Proposition}[section]
\newtheorem{theorem}[proposition]{Theorem}
\newtheorem{definition}[proposition]{Definition}
\newtheorem{lemma}[proposition]{Lemma}

\newtheorem{corollary}[proposition]{Corollary}
\newtheorem{remark}[proposition]{Remark}

\def\>#1{{\bf #1}}             

\def\dim{\hbox{{\rm dim }}}
\def\codim{\hbox{{\rm codim }}}
\def\Gr{\hbox{{\rm Gr}}}
\def\GL{\hbox{{\rm GL}}}
\def\Bi{\hbox{{\rm Bi}}}
\def\rank{\hbox{{\rm rank}}}
\def\vol{\hbox{{\rm vol}}}
\begin{document}

{\mbox{   }}

\vskip 2cm

\centerline{\Large\bf ALMOST HOLOMORPHIC EMBEDDINGS}

\bigskip

\centerline{\Large\bf IN GRASSMANNIANS  WITH APPLICATIONS}

\bigskip

\centerline{\Large\bf TO SINGULAR SYMPLECTIC SUBMANIFOLDS.}

\medskip

\vskip 2cm

\centerline{\sc V. Mu\~noz*, F. Presas** and I. Sols**}

\vskip 1.5cm

{\parindent 0cm{\it * Departamento de Matem\'aticas, 
Universidad Aut\'onoma de Madrid, 28049, Madrid, Spain.}}

{\parindent 0cm{\it ** Departamento de \'Algebra, Universidad 
Complutense de Madrid, 28040 Madrid, Spain.}}

\vskip 2cm

\centerline{\sc Abstract}
{\small In this paper we use Donaldson's approximately holomorphic 
techniques to build embeddings of a closed symplectic manifold with
symplectic form of integer class in
the grassmannians $\Gr(r,N)$. We assure that these embeddings are 
asymptotically holomorphic in a precise sense. We study first the 
particular case of $\CP^N$ obtaining control on $N$ and by a simple 
corollary we improve in a sense a classical result about symplectic 
embeddings \cite{Ti77}. The main reason of our study is the 
construction of singular determinantal submanifolds as the intersection 
of the embedding with certain ``generalized Schur cycles'' defined 
on a product of grassmannians. It is shown that
the symplectic type of these submanifolds is quite more general that
the ones obtained by Auroux \cite{Au97} as zeroes of ``very ample'' 
vector bundles.}

\vskip 1.5cm

\newpage

\section{Introduction and statement of the main results}
\label{introduction}

 Let $(M,\omega)$ be a symplectic manifold of integer class, 
 i.e.\ $[\omega/2\pi] \in H^2(M;\R)$ lifts to an integer 
 cohomology class. Such symplectic manifold has an associated 
 line bundle $L$ with first Chern class $c_1(L)=[\omega/2\pi]$,
 which is equipped with a connection $\nabla$ of curvature 
 $-i\omega$.

 In his outbreaking work \cite{Do96} S. Donaldson proved the 
 existence of symplectic submanifolds of $M$ that realize the 
 Poincar\'e dual of a large enough integer multiple of 
 $[\omega/2\pi]$. These are constructed as zero sets of 
 appropriate sections of $L^{\ox k}$.
 This extends a classical result in K\"ahler geometry 
 saying that $L$ is ample, so $L^{\otimes k}$ has holomorphic sections
 with smooth holomorphic, and so symplectic, zero sets.

 Later on, D. Auroux and R. Paoletti have proved independently an extension of 
 Donaldson's theorem, where now more symplectic submanifolds are constructed 
 as the zero sets of asymptotically holomorphic sections of vector 
 bundles.  These bundles are obtained by tensoring an arbitrary 
 complex bundle with large powers of the canonical line bundle 
 $L$ \cite{Au97}, \cite{Au99}, \cite{Pa99}.  In his paper,
 D. Auroux also shows that, asymptotically, all the sequences of 
 submanifolds constructed from a given
 vector bundle $E$ are isotopic. 
 (For a summary of these results see for example the review 
 paper \cite{Do98}.) 

 The key idea to understand these works is the concept of
 ampleness of a complex holomorphic bundle. This concept allows 
 the flexibilization of the bundles in the holomorphic category 
 by means of increasing their curvatures. Donaldson
 \cite{Do96} has translated the definition of ampleness to the symplectic
 category. For this he studies the asymptotical behaviour of sequences of
 sections of the bundles $L^{\ox k}$. Similarly, 
 the important point in our work is the definition of the concept of 
 asymptotic holomorphicity for sequences of embeddings constructed
 from very ample linear systems defined
 over vector bundles more and more twisted.

 The change to the non-integrable setting is controlled by this 
 concept. To define it we need to fix 
 a compatible almost complex structure $J$ in $(M,\omega)$. So the 
 pair $(\omega, J)$ gives a metric $g$ in the tangent bundle. We 
 have a sequence of metrics $g_k=kg$ indexed by integers $k\geq 1$.

\begin{definition} \label{def:ah}
 Let $X$ be a Hogde manifold with complex structure $J_0$. 
 Let $\g>0$. A sequence of embeddings $\phi_k:M\to X$ is 
 $\gamma$-asymptotically holomorphic 
 if it verifies the following conditions:
\begin{enumerate}
  \item $d\phi_k:T_xM \to T_{\p_k(x)}X$ has a left 
    inverse $\theta_k$ of norm less than $\gamma^{-1}$ at every 
    point $x\in M$. (The norm is taken with respect to the metric $g_k$.)
  \item $|(\phi_k)_*J- J_0|_{g_k}= O(k^{-1/2})$ on the subspace 
    $(\p_k)_* T_xM$.
  \item $|\nabla^p \phi_k|_{g_k}= O(1)$ and 
    $|\nabla^{p-1} \bbd \phi_k|_{g_k}= O(k^{-1/2})$, 
    for all $p\geq 1$.
\end{enumerate}
 A sequence of embeddings is asymptotically holomorphic if there is
 some $\g>0$ such that it is $\g$-asymptotically holomorphic.
\end{definition}

 The first important result is a generalization to the symplectic category of
 the classical Kodaira's embedding Theorem:

\begin{theorem} \label{existence0}
 Given $(M, \omega)$ a closed symplectic $2n$-dimensional manifold of 
 integer class endowed with a compatible almost complex structure, there
 exists an asymptotically holomorphic sequence of embeddings 
 $\phi_k:M \to \CP^{2n+1}$ with $\p^*_k[\o_{FS}]=[k\o]$. Moreover, given 
 two such sequences of embeddings asymptotically holomorphic with respect
 to two compatible almost complex structures, then 
 they are isotopic for $k$ large enough.
\end{theorem}

A sharper, in a sense, result than this has been obtained independently by 
Bortwick and
Uribe in \cite{BU99} using completely different ideas. Their result
also obtains control in the symplectic part (equivalently in the metric 
part) allowing to obtain asymptotically holomorphic embeddings which are also
asymptotically symplectic. Their approach is based on ideas coming
from Tian \cite{Ti90} to solve the problem in the K\"ahler case.

Our main interest for proving Theorem \ref{existence0} is given by the
possibility of studying ``projective symplectic geometry''. We mean by 
this the study of sequences of asymptotically holomorphic submanifolds, 
namely obtained as images of asymptotically holomorphic embeddings, 
in the projective space. The strength of this approach is shown in the
following

\begin{theorem} \label{good_inter0}
 Let $\phi_k$ be an asymptotically holomorphic sequence of embeddings 
 in $\CP^{2n+1}$ with $\p^*_k[\o_{FS}]=[k\o]$ and let $\epsilon>0$.
 Let us fix a holomorphic submanifold $N$ in $\CP^{2n+1}$. Then 
 there exists an asymptotically holomorphic sequence of
 embeddings $\hat{\phi}_k$, at distance at most $\epsilon$ in $C^r$-norm 
 from the initial sequence and verifying that 
 $\hat{\phi}_k(M)\bigcap N$ is symplectic for $k$ large enough.
\end{theorem}

With the notations introduced in Section \ref{section2} we will precise a 
little more the precedent result, assuring that $M \bigcap 
\hat{\phi}_k^{-1}(N)$ is a sequence of 
asymptotically holomorphic submanifolds.

We will see that this result will imply a projective version of the symplectic
Bertini's Theorem proved in \cite{Do96}. But the constructive method could 
allow to find more general types of symplectic submanifolds. This is shown
in a more general situation. For this we generalize Theorem \ref{existence0} 
to the grassmannian case.

\begin{theorem} \label{existence20}
 Let $(M, \omega)$ be a closed symplectic $2n$-dimensional manifold of integer
 class endowed with a compatible almost complex structure.
 Suppose also that we have a rank $r$ hermitian vector bundle with connection,
 and that $N>n+r-1$ and $r(N-r)>2n$. Then there exist an asymptotically 
 holomorphic sequence of embeddings $\phi_k:M \to \Gr(r, N)$ with $\p^*_k
 \SU=E\ox L^{\ox k}$, where $\SU\to \Gr(r,N)$ is the universal rank $r$ 
 bundle over the grassmannian. Moreover, given two such sequences of
 embeddings asymptotically holomorphic with respect to two compatible 
 almost complex structures, then they are isotopic for $k$ large enough.
\end{theorem}

In Section \ref{determinantal} we will take profit of this result to 
extend the construction of determinantal submanifolds to the symplectic
category in the following way.

\begin{definition} \label{defn:determ}
 Let $M$ be a differentiable manifold and let $E,F$ be complex vector bundles 
 over $M$. Given a morphism of vector bundles $\varphi: E\to F$, the 
 $r$-determinantal set $\Sigma^r(\varphi)$ is defined as
$$
 \Sigma^r(\varphi)= \{ x \in M \big| \rank\, \varphi_x=r \}. 
$$
\end{definition}

In the smooth category we can find
for any morphism $\varphi: E\to F$, 
another morphism $\hat{\varphi}: E\to F$ arbitrarily close to $\varphi$ 
in $C^p$-norm, such that
$\Sigma^r(\hat{\varphi})$ is a smooth submanifold in $M$
of codimension $2(r_e-r)(r_f-r)$, where $r_e$ and $r_f$ are the ranks of
$E$ and $F$, respectively (if this number is greater than the dimension of $M$
then the set is empty). 
There exists a similar result in the algebraic category if the vector bundle 
$E^* \otimes F$ is very ample. 
Our objective will be to adapt the algebraic discussion to the 
symplectic category to prove

\begin{theorem} \label{main_deter}
 Let $(M, \omega)$ be a closed symplectic manifold of integer class. 
 Let $E$ and $F$ be hermitian vector bundles of rank $r_e$ and $r_f$,
 respectively. Then, for $k$ large enough, there exists a morphism
 $\varphi_k:E\otimes (L^*)^{\otimes k}
 \to F\otimes L^{\otimes k}$ verifying that
\begin{enumerate}
 \item $\Sigma^r(\varphi_k)$ is an open symplectic submanifold of $M$.
 \item $\codim \Sigma^r(\varphi_k)= 2(r_e-r)(r_f-r)$. The set of manifolds 
  $\{\Sigma^r(\varphi_k)\}_r$
  constitutes a stratified submanifold, called determinantal submanifold.
\end{enumerate}
 Moreover, given two stratified determinantal submanifolds constructed 
 following the process described in the proof then there exists an ambient
 isotopy making the $r$-determinantal
 submanifolds associated to each stratified submanifold isotopic.
\end{theorem}

Theorem \ref{main_deter} was the original motivation of this paper.
The idea of studying this kind of submanifolds is inspired in algebraic 
geometry. Note that in algebraic geometry the manifolds constructed as 
zeroes of sections of vector bundles have many topological restrictions, 
namely they satisfy the Lefschetz hyperplane Theorem,
their Chern classes are very special, etc. So 
the set of submanifolds of a given manifold constructed in this way
is very special in the set of all the submanifolds.
However the determinantal submanifolds are very generic in the set
of submanifolds. For instance, every codimension $2$ submanifold of an
algebraic manifold can be constructed as the determinantal degeneration 
loci of certain bundle homomorphism \cite{Vo78}.

An obvious guess is that in symplectic geometry things are
similar. Recall that the most general submanifolds constructed using
asymptotically holomorphic techniques, prior
to Theorem \ref{main_deter} are the Auroux' ones \cite{Au97}. These 
are zeroes of sections of vector bundles, so its topological 
properties are very special.
In fact, Auroux cannot easily assure that these submanifolds are 
different from the ones constructed by Donaldson in \cite{Do96}. In 
Subsection \ref{chern} we compute some Chern numbers of determinantal
submanifolds showing that they are clearly different from the 
Chern numbers of Auroux' and Donaldson's submanifolds. So the symplectic 
type, and even the topological type, of the constructed submanifolds 
is necessarily different. This shows that the class of determinantal 
submanifolds is far more general.

Remark that, in any case, all the precedents results are obtained by means
of twisting vector bundles with large powers of the line bundle $L$. So
the submanifolds constructed in this way are quite special. It would be
desirable to avoid this restriction, but this generalization cannot
be made with the Donaldson's techniques developped in \cite{Do96}.

{}From a symplectic point of view determinantal submanifolds 
are also interesting. They constitute a step in the study
of singular symplectic submanifolds following the program sketched by Gromov 
\cite{Gr86}. Donaldson and Auroux have attacked this question in
\cite{Do99} and \cite{Au99}. Donaldson studies the local symplecticity of the 
fibers of asymptotically holomorpic applications $f:\C^n \to \C$ at a 
neigborhood of a critical point, it is solved by a local perturbation 
argument. The conclusion of Donaldson's work is that the topological 
behaviour of that kind of functions is similar to the holomorphic Morse 
functions. Auroux studies the local symplecticity of asymptotically 
holomorphic applications $f:\C^2 \to \C^2$ at the neighboorhood of a 
critical point, showing that are topologically equivalent to one of the 
two generic models of a holomorphic application \cite{Ar82}.
{}From this point of view Theorem \ref{main_deter} can be considered, in part,
an extension of these results to generic singularities. 

The organization of the paper is as follows. In Section \ref{section2}
we will give the basic ideas of the Donaldson-Auroux' theory needed
in our work and prove Theorem \ref{existence0}.
In Section \ref{projective} we prove Theorem \ref{good_inter0}. For this
we explain some euclidean notions concerning the estimation of
angles between subspaces. In Section \ref{emb_grass} we generalize all the
discussion to the case of the grassmannian embbedings, proving Theorem
\ref{existence20}. This allows us to prove Theorem \ref{main_deter} in Section 
\ref{determinantal} and to analyze the topological properties of the 
constructed submanifolds.

\vskip 5mm

\noindent {\bf Acknowledgements. }
We want to acknowledge D. Auroux, S. Donaldson and R. Paoletti
for his kindness communicating us their results.  Also we thank
A. Ibort and D. Mart\'{\i}nez for a lot of interesting discussions.
Second named author was conducting his research financed by the 
Ph.D.\ program of the Consejer\'{\i}a de Educaci\'on de Madrid.

\section{Asymptotically holomorphic embeddings in projective space}
\label{section2}
As in the introduction, let $(M,\omega)$ be a symplectic
manifold of integer class with associated line bundle $L$
and a compatible almost complex structure $J$.
In the K\"ahler setting this line 
bundle supports a holomorphic structure and it is ample in the 
algebraic geometry sense, i.e.\ $L^{\otimes k}$ has a lot of 
holomorphic sections. This allows to embed $M$ in the 
projective space $\CP^N$, for some $N$. 
In this Section we shall extend this classical result to the
symplectic case inspired in the ideas of \cite{Do96}, thus
proving Theorem \ref{existence0}.

\subsection{Asymptotically holomorphic sequences}
In this Subsection we collect the relevant results of the 
asymptotically holomorphic theory, as stated by D. Auroux 
in \cite{Au99}, that we shall use extensively along this work.

\begin{definition}
A sequence of sections $s_k$ of hermitian bundles $E_k$ 
with connections on $M$ is 
called asymptotically $J$-holomorphic if there exist constants 
$(C_p)_{p\in \N}$ such that, for all $k$, at every point of $M$, 
$|s_k|\leq C_0$, $|\nabla^p s_k|\leq C_p$ for all $p\geq 0$, 
and $|\nabla^{p-1} \bar{\partial} s_k|\leq C_p k^{-1/2}$ for all 
$p\geq 1$. The norms are evaluated with respect to the metrics $g_k$.
\end{definition}

In Donaldson's first work on the subject \cite{Do96}, 
$E_k=L^{\otimes k}$. In that work Donaldson imposed an
additional condition of improved transversality to the sequence of 
sections to assure that its zero sets are symplectic submanifolds 
for $k$ large enough. This condition is stated as follows.

\begin{definition} \label{skd-trans}
A section $s_k$ of the line bundle $L^{\otimes k}$ is said to be 
$\eta$-transverse to $0$ if for every point $x\in M$ such that 
$|s_k(x)|<\eta$ then $|\nabla s_k(x)|>\eta$.
\end{definition}

If we get an asymptotically $J$-holomorphic sequence $s_k$ of 
sections of $L^{\otimes k}$ verifying that all of them are 
$\eta$-transverse to $0$, with $\eta>0$ independent of $k$ 
then we can assure that $|\partial s_k(x)|>|\bar{\partial} s_k(x)|$ 
if $x$ is a zero of $s_k$, for $k$ large enough. A simple 
linear algebra argument assures that the zeroes of $s_k$ 
are symplectic submanifolds for $k$ large enough.

In \cite{Au97} D. Auroux extended the notion of transversality to
the case of higher rank bundles. Let $E$ be a rank $r$ hermitian
bundle with connection.

\begin{definition}
 A section $s_k$ of the bundle $E\ox L^{\ox k}$ is $\eta$-transverse
 to $0$ if for every $x\in M$ such that  
 $|s_k(x)|<\eta$ then $\nabla s_k(x)$ has
 a right inverse $\theta_k$ such that $|\theta_k|<\eta^{-1}$.
\end{definition}

We name {\it universal constant} to a number which only depends on the 
manifold geometry and on the constants involved in the data given
to start with, i.e.\ a number independent of $k$ and the point 
$x\in M$. Similarly a {\it universal polynomial} is a polynomial
only depending on the geometry of the manifold
and on the constants provided in the original data.
Donaldson uses highly localized asymptotically holomorphic sections, 
verifying the following definition.

\begin{definition}
A sequence of sections $s_k$ of hermitian bundles $E_k$ with
connections has Gaussian decay in 
$C^r$-norm away from the point $x\in M$ if there exists a universal 
polynomial $P$ and a universal constant $\lambda>0$ such that for 
all $y\in M$, $|s(y)|$, $|\nabla s(y)|_{g_k}$, $\ldots$, 
$|\nabla^r s(y)|_{g_k}$ are bounded by $P(d_k(x,y))\exp 
(-\lambda d_k(x,y))$. Here $d_k$ is the distance associated 
to the metric $g_k$.
\end{definition}

The starting point for Donaldson's construction is the following existence 
Lemma.

\begin{lemma} [\cite{Do96, Au97}] \label{localized}
 Given any point $x\in M$, for $k$ large enough, there exist asymptotically 
 holomorphic sections $\sref$ of $L^{\otimes k}$ over $M$ 
 satisfying the following bounds: 
 $|\sref|>c_s$ at every point of a ball of $g_k$-radius
 $1$ centered at $x$, for some universal constant $c_s>0$; 
 the sections $\sref$
 have Gaussian decay away from $x$ in $C^r$-norm.

 Moreover, given a one-parameter family of compatible almost-complex 
 structures $(J_t)_{t\in[0,1]}$, there exist one-parameter families 
 of sections $s_{t,k,x}^{\text{\scriptsize{ref}}}$
 which depend continuously on $t$ 
 and satisfy the same precedent properties. \hfill $\Box$
\end{lemma}

The proof of this Lemma uses in particular a refined version of 
Darboux' Theorem 
taking into account the holomorphic structure, which we also 
enunciate for later use.

\begin{lemma}[Lemma 3 in Chapter 3 of \cite{Au99}] \label{trivializa}
 Near any point $x\in M$, for any integer $k\geq 1$, there exist local 
 complex Darboux coordinates $(z_k^1, \ldots, z_k^n)=\Phi_k: (M,x) 
 \to (\C^n, 0)$ for the symplectic structure $k\omega$ such that the 
 followings bounds hold universally: $|\Phi_k(y)|^2=O(d_{k}(x,y)^2)$
 on a ball $B_{g_k}(x,c)$ of universal radius $c$ around $x$; 
 $|\nabla^r \Phi_k^{-1}|_{g_k}= O(1)$ for all $r\geq 1$ on a ball 
 $B(0,c')$ of universal radius $c'$ around $0$; and, with
 respect to the almost-complex structure $J$ on $X$ and the 
 canonical complex structure $J_0$ on $\C^n$, 
 $|\bar{\partial} \Phi_k^{-1}(z)|_{g_k}= O(k^{-1/2}|z|)$ and $|\nabla^r
 \bar{\partial} \Phi_k^{-1}|_{g_k}=O(k^{-1/2})$ for all $r\geq 1$ 
 on $B(0,c')$.

 Moreover, given a one-parameter continuous family of compatible 
 $(J_t)_{t\in [0,1]}$ and a continuous family of points 
 $(x_t)_{t\in[0,1]}$, there exists a continuous family of Darboux 
 coordinates $\Phi_{t,k}$ satisfying the same estimates and
 depending continuously on $t$.
\end{lemma}

{\bf Proof.}
 In \cite{Au99} the result is stated only for the case $n=2$ but the
 proof extends to the case $n>2$ trivially.
\hfill $\Box$

In \cite{Au99} D. Auroux used three asymptotically holomorphic sections to
set up a projection from a symplectic $4$-manifold $M$ to $\CP^2$. 
To control the behaviour of this projection he needs to assure global 
transversality conditions between the sections. He developes
a very useful scheme to pass from local transversality conditions to 
global ones by means of a globalization process inspired in the results 
of \cite{Do96}. Now we explain his idea to formalize Donaldson's techniques.

\begin{definition} \label{CrP}
A family of properties $\CrP(\epsilon, x)_{x\in M, \epsilon>0}$ of 
sections of bundles over $M$ is local
and $C^r$-open if, given a section $s$ satisfying $\CrP(\epsilon,x)$,
any section $\sigma$ such that $|\sigma(x)-s(x)|_{C^r}< \eta$ satisfies 
$\CrP (\epsilon-C\eta ,x)$, where $C$ is universal.
\end{definition}

For example, the property $|s(x)|>\epsilon$ is local and $C^0$-open.
The property that $s$ be $\epsilon$-transverse to $0$ at a
point $x$ is local and $C^1$-open.

\begin{proposition}[Proposition 3 in Chapter 3 of \cite{Au99}] 
\label{globalizate}
 Let $\CrP(\epsilon,x)_{x\in M, \epsilon >0}$ be a local and 
 $C^r$-open family of properties of sections of vector bundles 
 $E_k$ over $M$. Assume that there exist universal constants 
 $c$, $c'$, $c''$ and $p$ such that given any $x\in M$, any small 
 $\delta>0$, and asymptotically holomorphic sections $s_k$ of $E_k$, 
 there exist, for all large enough $k$, asymptotically holomorphic 
 sections $\tau_{k,x}$ of $E_k$ with the following properties:
\begin{enumerate}
  \item $|\tau_{k,x}|_{C^r, g_k}<c''\delta$.
  \item The sections $\frac{1}{\delta}\tau_{k,x}$ have Gaussian 
    decay away from $x$ in $C^r$-norm.
  \item The sections $s_k+\tau_{k,x}$ satisfy the property 
    $\CrP(\eta, y)$ for all $y\in B_{g_k}(x,c)$, with $\eta= 
    c'\delta (\log( \delta^{-1}))^{-p}$.
\end{enumerate}
 Then, given any $\alpha>0$ and asymptotically holomorphic sections 
 $s_k$ of $E_k$, there exist, for all large enough $k$, asymptotically 
 holomorphic sections $\sigma_k$ of $E_k$ such that 
 $|s_k-\sigma_k|_{C^r, g_k}< \alpha$ and the sections $\sigma_k$ 
 satisfy $\CrP(\epsilon, x)$ for all $x\in M$ with $\f>0$ independent
 of $k$.

 Moreover, the result holds for one-parameter families of sections, 
 provided the existence of sections $\tau_{t,k,x}$ satisfying 
 properties $1$, $2$ and $3$ and depending continuously on $t$.
\end{proposition}

{\bf Proof.}
 We only have added the constant $c''$ to the original statement in
 Proposition 3 in Chapter 3 of \cite{Au99}, which can be absorbed
 into the formula for $\eta$ just by enlarging $p$ universally.
\hfill $\Box$

The heart of these techniques is a series of local transversality
results which allow to apply Proposition \ref{globalizate}. These
results are based on ideas of complexity of real polynomials
coming from the real algebraic geometry. The most powerful result
is the following, proved in \cite{Do99}.

\begin{definition}
 A function $f:\C^n \to \C^r$ is $\s$-transverse to $0$ at a point 
 $x\in \C^n$ if it verifies at least one the following properties:
\begin{enumerate}
 \item $|f(x)|> \sigma$.
 \item $df(x)$ has a right inverse $\theta$ such that $|\theta|<\sigma^{-1}$.
\end{enumerate}
\end{definition}

\begin{proposition} {\rm (Theorem 12 in \cite{Do99})} \label{local_trans}
 There exists a universal integer $p$ verifying the following property: 
 for $0<\delta<\frac{1}{2}$ let $\sigma=\delta(\log (\delta^{-1}))^{-p}$.
 Let $f$ be a function with values in $\C^r$ defined over the ball $B^+=
 B(0,\frac{11}{10}) \subset \C^n$ satisfying the following bounds over $B^+$,
$$ 
  |f|\leq 1,~~~~ |\bbd f|\leq \s,~~~~ |\nabla \bbd f|\leq \s. 
$$
 Then there exists $w \in \C$ with $|w|<\d$ such that $f-w$ is 
 $\s$-tranverse to $0$ over the unit ball in $\C^n$. The same result holds 
 for one-parameter families of
 functions $f_t$ depending continously on $t\in [0,1]$, where
 we obtain a continuous path $w:[0,1] \to B(0, \delta)$. \hfill $\Box$
\end{proposition}

This Proposition is a generalization of Theorem 20 of \cite{Do96}, where
the case $r=1$ is proved. Later on D. Auroux in \cite{Au97,Au99}
extended the result to the parametric case with $r=1$ and to the case $r>m$
respectively. Proposition \ref{local_trans} covers all the range of 
possibilities. We mention also that in \cite{IMP99} the result is
refined to control the derivatives of the path $w_t$ allowing so
a generalization to the contact case of the asymptotically holomorphic
techniques.

\subsection{Asymptotically holomorphic embeddings in $\CP^{2n+1}$}
Through this Section we will study the existence of 
asymptotically holomorphic embeddings of a
closed symplectic manifold $(M,\o)$ of integer class and
dimension $2n$, endowed with a compatible almost complex structure $J$,
in the projective space $\CP^{2n+1}$. 
In Section \ref{emb_grass} we will develop the techniques to 
study the more general grassmannian embeddings. We want to prove the 
following

\begin{theorem} \label{existence1}
 Given an asymptotically $J$-holomorphic sequence of sections $s_k$ 
 of the vector bundles $\C^{2n+2}\otimes L^{\otimes k}$ and $\alpha>0$ 
 then there exists another sequence $\sigma_k$ verifying that:
\begin{enumerate}
  \item $|s_k-\sigma_k|_{C^1,g_k}< \alpha$.
  \item $\PP(\s_k)$ is an asymptotically holomorphic sequence of 
     embeddings in $\CP^{2n+1}$ for $k$ large enough.
  \item $[k\omega]=[\phi_k^* \omega_{FS}]$.
\end{enumerate}
 Moreover, let us have two asymptotically holomorphic sequences 
 $\phi_k^0$ and $\phi_k^1$ of embeddings in $\CP^{2n+1}$, with 
 respect to two compatible almost complex structures. Then for
 $k$ large enough, there exists an isotopy of asymptotically holomorphic 
 embeddings $\phi_k^t$ connecting $\phi_k^0$ and $\phi_k^1$. 
\end{theorem}

This result gives a proof of Theorem \ref{existence0}.
We shall proceed by steps to obtain asymptotically holomorphic 
embeddings of $M$ into $\CP^{2n+1}$.

\begin{definition}
 A sequence of asymptotically $J$-holomorphic sections $s_k$ of the 
 vector bundles $\C^{2n+2} \otimes L^{\otimes k}$ is $\gamma$-projectizable 
 if for all $x\in M$, $|s_k(x)|>\gamma$.
\end{definition}

This is a sufficient condition to get a map to $\CP^{2n+1}$ defined 
as $\phi_k=\PP (s_k): M \to \CP^{2n+1}$, as the $\gamma$-projectizability
assures that the sections $s_k=(s_k^0, \ldots, s_k^{2n+1})$ are not 
simultaneously zero and so the $\PP$ operator is well defined. To get
local injectivity we need to impose the following.

\begin{definition}
 Let $s_k$ be a sequence of asymptotically $J$-holomorphic 
 $\gamma$-projectizable sections of the vector bundles 
 $\C^{2n+2}\otimes L^{\otimes k}$ for some $\g>0$ and let 
 $0\leq l\leq n$. Then $s_k$ is 
 $\eta$-generic of order $l$, with $\eta>0$, if
 $|\bigwedge^l \partial \PP(s_k)(x)|_{g_k}>\eta$ for all $x\in M$.
 For $l=0$ the condition is vacuus.
\end{definition}

We have the following result that will be proved in the
following two Subsections.

\begin{proposition} \label{key}
 Let $s_k$ be an asymptotically $J$-holomorphic sequence of sections of
 the vector bundles $\C^{2n+2}\otimes L^{\otimes k}$ and $\a>0$. Then 
 there exists another asymptotically holomorphic sequence $\s_k$ verifying:
\begin{enumerate}
  \item $|s_{k}-\sigma_{k}|_{C^1,g_k}< \alpha$.
  \item $\sigma_{k}$ is $\gamma$-projectizable and $\gamma$-generic 
    of order $n$ for some $\gamma>0$.
\end{enumerate}
 Moreover, the result holds for one-parameter families of sections
 where the sections and almost complex structures depend continuously 
 on $t\in [0,1]$.
\end{proposition}

With this result we can give the proof of Theorem \ref{existence1}.

{\bf Proof of Theorem \ref{existence1}.}
We first prove the existence result. The last property is obvious since 
the hyperplane bundle of $\CP^{2n+1}$ restricts by construction to 
$L^{\otimes k}$. Let us begin with an asymptotically $J$-holomorphic 
sequence $\s_k$ of sections of the bundles $\C^{2n+2}\otimes 
L^{\otimes k}$. Now we perturb it using Proposition \ref{key} 
to obtain an asymptotically holomorphic sequence $s_k$ with
$|s_{k}-\sigma_{k}|_{C^1,g_k}< \alpha$, which is $\g$-projectizable 
and $\g$-generic of order $n$, for some $\g>0$. We have only to 
check that the sequence $\phi_k=\PP (s_k)$ satisfies the required
properties in Definition \ref{def:ah}.
More specifically, we shall check that $\p_k$ is an 
immersion of $M$ in $\CP^{2n+1}$, for $k$ large. 
To get rid of the possible self-intersection we take into account 
that $2\, \dim M< \dim \CP^{2n+1}$ so we can make a generic 
$C^r$-perturbation of norm less than $O(k^{-1/2})$ to get an embedding 
keeping the asymptotic holomorphicity and the genericity of order 
$n$.

Choose a point $x\in M$. By a rotation with an element of $U(2n+2)$
acting on $\C^{2n+2}$, we can assure that $s_k(x)=(s_k^0(x), 
\ldots, s_k^{2n+1}(x))=(s_k^0(x),0,\ldots, 0)$.
The transformation is constant on $M$ and only produces a global 
isometric transformation of $\phi_k(M)$ in $\CP^{2n+1}$. 
Now using the $\gamma$-projectizable property we know that 
$|s_k^0(x)|\geq \gamma$. By the asymptotically holomorphic 
bounds of $s_k^0$ there is a universal $c$ such that
$|s_k^0|\geq \gamma/2$ on $B_{g_k}(x,c)$ for all $k$.
We define the application:
\begin{eqnarray*}
  f_k:  B_{g_k}(x,c) & \to & \C^{2n+1} \\
  y & \mapsto & \left( \frac{s_k^1(y)}{s_k^0(y)}, \ldots, 
   \frac{s_k^{2n+1}(y)}{s_k^0(y)}\right).
\end{eqnarray*}
This application can be written as $f_k= \Phi_0 \circ \phi_k$, 
where $\Phi_0$ is the 
standard trivialization application in $\CP^{2n+1}$ defined for the chart 
$U_0=\{ x=[x_0, \ldots, x_{2n+1}]\big| x_0\neq 0\}$. It is well known 
that $\Phi_0$ is an isometry at the point $[1,0, \ldots,0]$ if we use 
the standard metric structure of $\C^{2n+1}$. So we can compute the bounds 
required in Definition \ref{def:ah} using $f_k$ instead of $\phi_k$. 
The asymptotic holomorphicity of $s_k$ and the bound 
$|s_k^0|\geq \gamma/2$ imply that $|\nabla^p f_k(x)|=O(1)$ 
and $|\nabla^p \bar{\partial} f_k(x)|=O(k^{-1/2})$, for $p \geq 0$.
This proves condition 3 in Definition \ref{def:ah}. 

Now we pass to the issue of the existence of a left inverse. We have the
decomposition
$$ 
  \bigwedge\nolimits^n d\phi_k= \bigwedge\nolimits^n \partial \phi_k 
  + O(k^{-1/2}), 
$$
where the last term is obtained thanks to $|\bbd\p_k|_{g_k}=O(k^{-1/2})$.
By the $\g$-genericity of order $n$ of $\p_k$, 
$|\bigwedge^n \bd\phi_k|_{g_k} \geq \g$, so
$|\bigwedge^n d\phi_k|_{g_k} \geq \g/2$ for $k$ large. Let
$$
  \hat{\theta}_k = (d\p_k)^{-1}: (\p_k)_* T_xM \to  T_xM.
$$
By the asymptotic holomorphicity condition,
we have $|d\p_k|_{g_k} \leq C_0$ for a universal constant $C_0$, so
$|\hat{\theta}_k| \leq C \g^{-1}$ for another universal constant 
$C$. Now define $\theta_k=\hat{\theta}_k 
\circ \pr^{\perp}$, where $\pr^{\perp}$ is the orthogonal projection 
of $T_{\phi_k(x)}\CP^{2n+1}$ onto $(\phi_k)_* T_xM$ 
to get the sought right inverse (reducing $\g$ conveniently). 

Finally we compute the norm of $(\p_k)_* J-J_0:
(\p_k)_* T_xM \to T_{\p_k(x)} \CP^{2n+1}$. 
The expression can be written as
$$
  (\phi_k)_* J- J_0=d\phi_k J\hat{\theta}_k - J_0=
  (d\p_k+J_0 d\p_k J) J \hat{\theta}_k= 2\bbd \p_k J \hat{\theta}_k=O(k^{-1/2}),
$$
proving condition 2 in Definition \ref{def:ah}.

For the isotopy result we follow the ideas of \cite{Au97}. We need the 
following auxiliary result, which we prove in Subsection \ref{lifting}.

\begin{lemma} \label{lem:emb->sec}
  Let $\p_k:M\to \CP^{2n+1}$ be a sequence of asymptotically 
  holomorphic embeddings with $\p_k^*[\o_{FS}]=[k\o]$. 
  Then there exists a sequence of asymptotically holomorphic 
  sections $s_k$ of $\C^{2n+2}\ox L^{\ox k}$, for $k$ large enough, 
  which is $\g$-projectizable and $\g$-generic of order $n$, for some
  $\g>0$, such that $\p_k=\PP(s_k)$. The same holds for continuous
  one-parameter families of embeddings and compatible almost complex
  structures. 
\end{lemma}
 
Using Lemma \ref{lem:emb->sec}, we can suppose that $\p_k^i=
\PP(s_k^i)$, $i=0,1$, where $s_k^0$ and $s_k^1$ are 
two asymptotically holomorphic sequences which
are $\g$-projectizable and $\g$-generic of order $n$, $\g>0$. We construct 
the following family of sequences of asymptotically holomorphic 
sections:
$$
  s_k^t= \left\{ \begin{array}{lll} (1-3t)s_k^0, \qquad &  
       {\rm with} \ J_t=J_0, & t\in [0,1/3] \\
   0, & {\rm with} \ J_t=\text{Path}(J_0,J_1), \qquad & t\in [1/3,2/3] \\ 
   (3t-2)s_k^1, & {\rm with} \ J_t=J_1, & t\in [2/3,1]. \end{array} \right.
$$
Choose $\a>0$ such that any perturbation of $s_k^t$ of $C^1$-norm 
less than $\a$ is still $\g/2$-projectizable and $\g/2$-generic of
order $n$.
Applying Proposition \ref{key} to $s_k^t$ with this $\a$, we obtain a
family $\sigma^t_k$ which is $\eta$-projectizable and $\eta$-generic of
order $n$ for some $\eta>0$. We define the family of sequences of 
asymptotically holomorphic sections:
$$ 
  \tau_k^t= \left\{ \begin{array}{ll} 
   (1-3t)s_k^0+3t\sigma_k^0, & t\in [0,1/3] \\
   \sigma_k^{3t-1}, & t\in [1/3,2/3] \\ 
   (3t-2)s_k^1+(3-3t)\sigma_k^1, \qquad & t\in [2/3,1]. \end{array} \right. 
$$
These are $\f$-projectizable and $\f$-generic of order $n$ sequences of
sections, with $\f=\min\{\g/2,\eta\}$, so that $\p^t_k=\PP(\tau_k^t)$ are 
asymptotically holomorphic embeddings
(maybe after a further small perturbation to get rid of self-intersections). 
This implies that $\phi_k^0$ and $\phi_k^1$ are isotopic for
$k$ large enough.
\hfill $\Box$

An important corollary is the existence of symplectic embeddings 
of $M$. The following result is similar to \cite{Ti77}, but we do 
not obtain an exact symplectic embedding. On the other hand the 
dimension of the projective space is controlled in our case.

\begin{corollary} \label{symp_emb}
 Let $(M,\omega)$ be a closed symplectic manifold
 of dimension $2n$ with symplectic form of integer class. Then there 
exists a symplectic embedding 
 $\phi: M \to \CP^{2n+1}$ verifying that $k\omega= \phi^* \omega_{FS}$,
 for $k$ large enough.
\end{corollary}

{\bf Proof.}
 Take a $\g$-asymptotically holomorphic sequence $\phi_k$ of embeddings 
 of $M$ in $\CP^{2n+1}$. The key idea is that the linear segment of 
 forms $\omega_t$ joining two symplectic forms compatible with a fixed 
 $J$ is symplectic for every $t$. In our case we have this condition
 asymptotically. Define the family of $2$-forms in $M$ given by 
 $\omega_t=(1-t)k\omega+ t\phi_k^*(\omega_{FS})$, where $t\in [0,1]$. 
 All of them are cohomologous, so to apply Moser's trick \cite{MS94}
 we only need to prove that they are symplectic. Suppose that
 there exists $t$ such that $\omega_t$ is not symplectic. Then there 
 is a unitary tangent vector $v\in T_xM$, for some $x\in M$, such that 
 $\omega_t(v,w)=0$, for all $w\in T_xM$. In particular 
 $\omega_t(v,Jv)=0$. Now expanding this expression we obtain:
\begin{eqnarray*}
  \omega_t(v,Jv) & = & (1-t)k\omega(v,Jv)+t\phi_k^* \omega_{FS}(v, Jv) \\
   & = & (1-t)kg(v,v)+t\o_{FS}(d\p_k v, J_0 \bd\p_k v -J_0 \bbd \p_k v) \\
  & = & (1-t)kg(v,v)+tg_{FS}(d\p_k v,\bd\p_k v)-tg_{FS}(d\p_k v, \bbd
     \phi_k v) \\
  & = & (1-t)kg(v,v)+tg_{FS}(d\p_k v,d\p_k v)
    -2tg_{FS}(d\p_k v, \bbd \phi_k v) \\
  & = & (1-t)kg(v,v)+tg_{FS}(d\p_k v,d\p_k v)-tO(k^{-1/2}). 
\end{eqnarray*}
 Thanks to the $\g$-asymptotically holomorphic embeddings, we have that
 $g_{FS}(d\phi_k v, d \phi_k v) \geq \g^2$. So for $k$ large enough
 we get a contradiction. \hfill $\Box$

\subsection{Construction of $\gamma$-projectizable sections.}
Our objective is to prove the following perturbation result.

\begin{proposition}
 Let $s_k$ be an asymptotically $J$-holomorphic sequence of sections 
 of vector bundles $\C^{2n+2}\otimes L^{\otimes k}$. Then given 
 $\alpha>0$, there exists an asymptotically
 $J$-holomorphic sequence of sections $\sigma_k$ verifying:
\begin{enumerate}
  \item $|s_k-\sigma_k|_{C^1,g_k}< \alpha$.
  \item $\sigma_k$ is $\eta$-projectizable for some $\eta>0$.
\end{enumerate}
 Moreover, the result can be extended to continuous one-parameter 
 families of asymptotically $J_t$-holomorphic sequences $s_{t,k}$ 
 obtaining continuous one-parameter families 
 of sections $\sigma_{t,k}$ verifying the two precedent conditions.
\end{proposition}

{\bf Proof.}
The result is a simple generalization of Proposition 1 in \cite{Au99} 
where the result for $4$-manifolds is proved. The high dimensional 
case can be treated with the same techniques. 

 We will proceed by using the globalization argument described in 
 Proposition \ref{globalizate}. First we deal with the non-parametric
 case. For this we define the local and $C^0$-open property 
 $\CrP (\epsilon,x)$ as $|s_k(x)| > \epsilon$. Let $\d>0$. 
 We only need to find 
 for a point $x\in M$ a section $\tau_{k,x}$ with Gaussian decay 
 away from $x$, assuring that $s_k+\tau_{k,x}$ verifies $\CrP(\eta,
 y)$ in a ball of universal $g_k$-radius $c$, with $\eta=c' \d
 (\log (\d^{-1}))^{-p}$, $c'$ and $p$ universal constants.

 For this choose a section $\sref$ verifying the conditions 
 of Lemma \ref{localized}. Then we select $c=1$ (obviously, universal). 
 The lower bound of $\sref$ in the ball $B_x=B_{g_k}(x,1)$ let
 us define the application 
$$  
 f_{k,x}=\frac{s_k}{\sref}: B_x \to \C^{2n+2}.
$$
 Using the lower bound of $\sref$ together with the asymptotic
 holomorphicity of $s_k$ is easy to show that
\begin{equation} \label{bounds} 
  |f_{k,x}|< C, ~~~~ |\bar{\partial} f_{k,x}|< Ck^{-1/2}, ~~~~ 
  |\nabla \bar{\partial} f_{k,x}|< Ck^{-1/2},
\end{equation}
 where $C$ is a universal constant. With the aid of Lemma \ref{trivializa}
 we can build $f_k=f_{k,x} \circ \Phi_k^{-1}$ defined on a fixed ball 
 $B(0,c')\subset \C^n$. Scaling the coordinates by a universal constant 
 $\frac{11}{10}(c')^{-1}$ we can suppose that $f_k$ is defined on $B^+$. 
 In this ball, the bounds (\ref{bounds}) yield
\begin{equation}\label{boundstriv}
  |f_k|< C_0, ~~~~ |\bar{\partial} f_k|< C_0k^{-1/2}, ~~~~ 
  |\nabla \bar{\partial} f_k|< C_0k^{-1/2},
\end{equation}
 where $C_0$ is a universal constant. The application 
 $g_k=\frac{1}{C_0}f_k$ is in the hypothesis of Proposition \ref{local_trans} 
 and then there exists, for $k$ large enough, a number $w_k\in 
 B(0, \d)$ such that $|g_k-w_k|>\s=\delta(\log (\delta^{-1}))^{-p}$. 
 Therefore $|f_k-C_0w_k|>C_0\s$ on $B$. Now define 
 $\tau_{k,x}=-C_0 w_k\otimes \sref$, so that 
 $|\tau_{k,x}|_{C^r,g_k} <c'' \d$, for some universal constant $c''$.
 Using the lower bound of $\sref$ we obtain that 
 $|s_k+\tau_{k,x}|\geq c'\delta (\log (\delta^{-1}))^{-p}$,
 with $c'$ and $p$ universal constants. 
 Then Proposition \ref{globalizate} applies
 and the proof is concluded in the non-parametric case.

 The globalization to the one-parameter case is trivial because all 
 the ingredients in the proof can be easily chosen in a continuous way. 
\hfill $\Box$

\subsection{Inductive construction of sections $\gamma$-generic of order $l$}
Now we study the problem of perturbing a $\g$-projectizable sequence 
of sections to achieve genericity of order $n$. We shall do this in steps.
The result to be proved is the following

\begin{proposition} \label{generic}
 Let $s_k$ be an asymptotically $J$-holomorphic sequence of sections 
 of the vector bundles $\C^{2n+2}\otimes L^{\otimes k}$ which is 
 $\g$-projectizable and $\g$-generic of order $l$. Then given $\alpha>0$, 
 there exists an asymptotically $J$-holomorphic sequence of sections 
 $\sigma_k$ verifying:
\begin{enumerate}
  \item $|s_k-\sigma_k|_{C^1,g_k}< \alpha$.
  \item $\sigma_k$ is $\eta$-generic of order $l+1$ for some $\eta>0$.
\end{enumerate}
 Moreover, this can be extended to continuous one-parameter families 
 of asymptotically $J_t$-holomorphic sequences $s_{t,k}$ obtaining 
 continuous one-parameter families of sections $\sigma_{t,k}$ verifying 
 conditions 1 and 2.
\end{proposition}

{\bf Proof.}
 We construct local $1$-forms to control the perturbations. For this 
 at a neighborhood of a point $x\in M$ we fix local complex Darboux 
 coordinates $(z_k^1, \ldots, z_k^n)$ using Lemma \ref{trivializa}.
 As in proof of Theorem \ref{existence1}, by applying a unitary transformation
 to $\C^{2n+2}$, we can suppose that $s_k(x)=(s_k^0(x), 0, \ldots, 0)$. 
 Also there exists a ball with center $x$ and universal $g_k$-radius 
 $c$ on which $|s_k^0| \geq\g/2$.  
 We define, following Auroux' notations \cite{Au99}, a local basis of
 asymptotically holomorphic $1$-forms:
$$
  \mu_k^j= \partial \left( \frac{z_k^j \sref}{s_k^0} \right),
$$
 where $\sref$ are given by Lemma \ref{localized}.
 They have Gaussian decay away from $x$ thanks to the behaviour 
 of $\sref$. At $x$ they form an orthonormal basis of $T^*_xM$. 
 We use the trivialization $\Phi_0$ to define the application
\begin{eqnarray}
 f_k:  B_{g_k}(x,c) & \to & \C^{2n+1} \label{entorno} \\
 y & \mapsto & \left(\frac{s_k^1(y)}{s_k^0(y)}, \ldots, 
   \frac{s_k^{2n+1}(y)}{s_k^0(y)}\right), \nonumber
\end{eqnarray}
 which is almost an isometry on $B_{g_k}(x,c)$. 
 
 The case $l=0$ without parameters is the easiest. We say that a section 
 $\g/2$-projectizable verifies $\CrP (\epsilon, x)$ if 
 $|\bd \phi_k(x)|>\epsilon$. This property is local and open
 in $C^1$-sense. We are going to apply Proposition \ref{globalizate} 
 to assure the existence of a $\eta$-generic of order $1$ sequence 
 of sections arbitrarily near the given $s_k$ in $C^1$-norm, for
 some $\eta>0$. For this let $0<\d<\g/2c''$, $c''$ a universal
 constant whose precise value will appear later. We have to build a local 
 perturbation $\tau_{k,x}$ with $|\tau_{k,x}|<c''\d$ and Gaussian 
 decay to achieve the property 
 $\CrP(\eta, y)$ in a neighborhood of $x$ of universal $g_k$ radius $c$, 
 with $\eta=c'\d(\log(\d^{-1}))^{-p}$. 
 (As we only perturb with sections of $C^0$-norm less than $\g/2$ we 
 can assure that all the sections still have the property 
 $\g/2$-projectizable.)

 Fixing $x\in M$, we have the applications $f_k$ of \eqref{entorno}.
 It is easy to check that there is a ball of universal radius $c_0$ where 
$$ 
 \bd f_k = (u_k^{11}\mu_k^1+u_k^{12}\mu_k^2+ \cdots + u_k^{1n}\mu_k^{n}, 
 \,\ldots\, , u_k^{2n+1,1}\mu_k^1+\cdots + u_k^{2n+1,n}\mu_k^n ), 
$$
 for some $u_k^{ij}$. Then we obtain an application $u_k:B_{g_k}(x,c_0) 
 \to \C^{n\x (2n+1)}$. Using a complex Darboux chart we can trivialize 
 $B_{g_k}(x, c_0)$ to obtain (scaling the coordinates by an appropriate
 universal constant $C$) an application $\hat{u}_k: B^+ \to 
 \C^{n\x (2n+1)}$ which is asymptotically holomorphic by construction.
 So we can apply Proposition \ref{local_trans} to get
 $w_k'\in \C^{n\x (2n+1)}$ such that $|\hat{u}_k-w_k'|>\eta=
 \delta (\log (\delta^{-1}))^{-p}$ on $B$, where $|w_k'|<\delta$. 
 Rescaling and passing to the manifold we obtain that 
 $|u_k-Cw_k'|>C \delta (\log (\delta^{-1}))^{-p}$.
 We denote $w_k=Cw_k'$ and define the section 
 $\tau_{k,x}=-(0, w_k^{11}z_k^1\sref+ 
 w_k^{12}z_k^2\sref+ \cdots + w_k^{1n}z_k^n \sref, 
 \ldots,  w_k^{2n+1,1}z_k^1\sref+ \cdots +
 w_k^{2n+1,n}z_k^n \sref)$ of $\C^{2n+2}\ox L^{\ox k}$.
 This section verifies the properties required in 
 Proposition \ref{globalizate}. 

 To check the one-parameter case we have only to get a continuous
 family of unitary transformations verifying that 
 $s_{t,k}(x)=(s_{t,k}^0(x), 0, \ldots, 0)$ for all $t\in [0,1]$.
 This is clearly possible because of the contractibility of $[0,1]$.

 Now we pass to the case $l>0$. We define the following property for 
 sections $s_k$ which are $\g/2$-projectizable and  
 $\g/2$-generic of order $l$. A section $s_k$ has the property
 $\CrP (\epsilon, x)$ if $|\bigwedge^{l+1} \partial \PP s_k(x)|>\epsilon$. 
 This property is local and open in $C^1$-sense. For applying
 Proposition \ref{globalizate} we need to build, for $0<\d<\g/2c''C$, 
 a local
 perturbation $\tau_{k,x}$ with $|\tau_{k,x}|<c''\d$ and Gaussian 
 decay with the property $\CrP(\eta, y)$ in a neighborhood of $x$ 
 of universal $g_k$ radius $c$, with $\eta=c'\d(\log(\d^{-1}))^{-p}$.
 (Here $C$ is the constant of the $C^1$-openness of $\CrP(\f,x)$ in 
 Definition \ref{CrP}.)
 We define $f_k$ as in \eqref{entorno}. Then it is easy to see that
 there exists a universal constant $c$ such that 
$$
 \frac{|\bigwedge^{l+1} \partial \PP(s_k)|}{|\bigwedge^{l+1} \partial f_k|}
 > 1/2
$$
 on $B_{g_k}(x,c)$. So we can do the computations for the applications
 $f_k$. By a unitary transformation in $U(2n+1)$ (on $\C^{2n+2}$ fixing 
 $(1,0,\ldots,0)$) and other in $U(n)$ (on the complex Darboux coordinate 
 chart) we can assure that
\begin{equation} 
  \bd f_k(x)= \left( \begin{array}{ccccccc}
  u_k^{11}(x) & 0 & \ldots & & & \ldots & 0 \\
  0 & u_k^{22}(x) & 0 & \ldots & & & 0 \\
  0 & \ldots & \ddots & 0 & \ldots & & 0 \\
  0 & \ldots &0 & u_k^{nn}(x) & 0 & \ldots & 0 \end{array} \right), 
  \label{matriz}
\end{equation}
 where $|u_k^{11}(x)\cdots u_k^{ll}(x)|>\g/C'$, $C'$ a universal constant. 
 Shrinking $c$ if necessary we can assure that 
 $|(\partial f_k^1\wedge \cdots \wedge \partial f_k^l)_{\mu_k^1\wedge \cdots 
 \wedge \mu_k^l}|> \g/2C'$ for all the points of the 
 ball $B_{g_k}(x,c)$, where we denote by $(\partial f_k^1\wedge \cdots \wedge 
 \partial f_k^l)_{\mu_k^1\wedge \cdots \wedge \mu_k^l}$ the component of 
 $\partial f_k^1\wedge \cdots \wedge \partial f_k^l$ in the direction of 
 $\mu_k^1\wedge \cdots \wedge \mu_k^l$. This $l$-form is an element of 
 the basis composed by the $l$-wedge products of the $1$-forms $\mu_k^1, 
 \ldots, \mu_k^n$. In matrix form we are denoting the order $l$ left 
 upper minor of $\partial f_k$. Now we construct the $(l+1)$-form
$$ 
  \theta_k(y)= (\partial f_k^1\wedge \cdots \wedge \partial 
  f_k^l)_{\mu_k^1\wedge \cdots \wedge \mu_k^l} \wedge \mu_k^{l+1}. 
$$
 We can suppose that $|\theta_k|>c_s \g$ with $c_s>0$ a universal 
 constant. We also consider the following family of $(l+1)$-forms
$$
  M_k^p= (\partial f_k^1\wedge \cdots \wedge \partial f_k^l \wedge 
  \partial f_k^p)_{\mu_k^1 \wedge \cdots \wedge \mu_k^l \wedge \mu_k^{l+1}}, 
  \quad l+1\leq p \leq 2n+1.
$$
 These forms are components of $\bigwedge^{l+1} \partial f_k$. If we 
 perturb so that the norm of $M_k=(M_k^{l+1}, \ldots , M_k^{2n+1})$ is 
 bigger than $\eta=c'\d (\log (\d^{-1}))^{-p}$ then we have
 finished because if $|M_k|>\eta$ then 
 $|\bigwedge^{l+1} \partial f_k|> C_0\eta$ where
 $C_0$ is again a universal constant (using that
 the basis $\{ \mu_k^{i_1} \wedge \cdots \wedge 
 \mu_k^{i_{l+1}} \}_{1\leq i_1 < \cdots < i_{l+1}\leq n}$ 
 is almost orthogonal on the ball $B_{g_k}(x,c)$, in fact orthogonal 
 at $x$). 

 We define the following sequence of asymptotically holomorphic
 applications,
$$
g_k=(g_k^{l+1}, \ldots, g_k^{2n+1})=\left(
 \frac{M_k^{l+1}}{\theta_k}, \ldots, \frac{M_k^{2n+1}}{\theta_k}\right).
$$
 So we obtain, scaling the coordinates by universal constants if necessary,
 $\hat g_k: B^+ \to \C^{2n+1-l}$ which is asymptotically holomorphic thanks to 
 the lower bound of $\theta_k$ and to the asymptotic holomorphicity of 
 $M_k$ and $\theta_k$. We have that $n<2n+1-l$ and so we 
 can find $|w_k|<\delta$ such that 
 $|g_k-w_k|>\eta= \delta (\log (\delta ^{-1}))^{-p}$.
 Thus we obtain that $|(M_k^{l+1}-w_k^{l+1}\theta_k,
 \ldots, M_k^{2n+1}-w_k^{2n+1}\theta_k)|>c_s\g \eta$. Recall that all
 the constants depend
 on $\g$ and the asymptotic holomorphicity constants of
 $s_k$, so they are independent  of $x$ and $k$.
 The perturbation $-(w_k^{l+1} \theta_k, \ldots, w_k^{2n+1} \theta_k)$
 is achieved by adding the section
 $\tau_{k,x}=-(0, \stackrel{(l)}{\cdots}, 0, w_k^{l+1}z_k^{l+1}\sref,
 \ldots, w_k^{2n+1}z_k^{l+1}\sref)$ to $s_k$. 
 This section verifies the Gaussian decay
 bounds required in Proposition \ref{globalizate}
 and $|\tau_{k,x}|_{C^1,g_k}<c''\d$ for some universal constant $c''$. 
 This completes the proof in the non-parametric case.

 Now we pass to the one-parameter case. With appropriate continuous 
 unitary transformations, we may 
 assume that $s_{t,k}(x)=(s_{t,k}^0(x), 0, \ldots, 0)$ and that $\bd 
 f_{t,k}(x)$ is written as in \eqref{matriz}. The interval $[0,1]$ may
 be split in a finite number of subintervals $[t_i,t_{i+1}]$ such that,
 for every $x\in M$ and for each of the subintervals, there is
 a fixed order $l$ minor of $\bd 
 f_{t,k}(x)$ with norm bigger than $\g/C'$, for every $t$ in the
 subinterval. This allows to find global small perturbations of $s_{t,k}$
 in every $[t_i,t_{i+1}]$. Reducing $\a$ and enlarging $C'$ we may 
 suppose that the same happens to any perturbation of the original 
 $s_{t,k}$ at $C^1$-distance at most $\a$.

 Now work as follows. For the first subinterval, 
 consider $s^1_{t,k}= s_{t,k}$, $t\in [0,t_1]$. We find 
 a perturbation $\s^1_{t,k}$, $t\in [0,t_1]$,
 such that $|s^1_{t,k}-\s^1_{t,k}|<\a/2$ and $\s^1_{t,k}$ is 
 $\eta_1$-generic of order $l+1$, for some $\eta_1>0$. 
 Set $\s_{t,k}=\s^1_{t,k}$ for $t\in [0,t_1]$.
 In the second subinterval, perturb $s_{t,k}^2= s_{t,k}^1 +
 (\s^1_{t_1,k}-s^1_{t_1,k})$, $t\in [t_1,t_2]$, to find $\s^2_{t,k}$
 satisfying $|s^2_{t,k}-\s^2_{t,k}|<\a/4$ and $\s^2_{t,k}$ is 
 $\eta_2$-generic of order $l+1$, for some $\eta_2>0$.
 To glue this perturbation with the previous one puts
$$
 \s_{t,k}= \left\{ \begin{array}{ll} s_{t,k}^2+\frac{t-t_1}{\f} 
   (\s_{t,k}^2 -s_{t,k}^2), \qquad & t\in [t_1,t_1+\f] \\
      \s_{t,k}^2,  & t\in [t_1+\f,t_2]. \end{array} \right.
$$
 Here $\f>0$ is chosen so small that 
 $|s_{t,k}^2-\s_{t_1,k}^1|_{C^1}<\rho/2$, for $t\in [t_1,t_1+\f]$, 
 and we require also that the perturbation satisfies
 $|s^2_{t,k}-\s^2_{t,k}|<\rho/2$, where $\rho>0$ is a number such that
 any perturbation of $\s_{t_1,k}^1$ of $C^1$-norm less than $\rho$ 
 is $\eta_1/2$-generic of order $l+1$. This defines $\s_{t,k}$ for
 $t\in [0,t_2]$ already.

 Proceeding in this way we finally find $\s_{t,k}$, $t\in [0,1]$, 
 which is $\eta$-generic of order $l+1$, for some $\eta>0$, with
 $|\s_{t,k}-s_{t,k}|< \a$.
\hfill $\Box$

\subsection{Lifting asymptotically holomorphic embeddings}
\label{lifting}
In this Subsection we aim to prove that the sequences of asymptotically
holomorphic embeddings into $\CP^{2n+1}$ that we are considering in
Theorem \ref{existence0} come always from asymptotically holomorphic
sequences of sections $s_k$ of $\C^{2n+2}\ox L^{\ox k}$ which are 
$\g$-projectizable and $\g$-generic of order $n$, for some $\g>0$
(at least for $k$ large).

{\bf Proof of Lemma \ref{lem:emb->sec}.}
Suppose that we have a sequence of $\g$-asymptotically holomorphic
embeddings $\p_k:M\to \CP^{2n+1}$, for some $\g>0$, with
$\p_k^*\SU=L^{\ox k}$. Here $\SU$ is the hyperplane line bundle
defined over the projective space. The dual of $\SU$ is the 
universal line bundle
$$
  \SE=\{(l,s) \big| s \in l \} \subset \CP^{2n+1}\x \C^{2n+2}=
  \underline{\C}^{2n+2},
$$
interpreted as a sub-bundle of the trivial bundle $\underline{\C}^{2n+2}$.

Consider the following sequence of line bundles, $E_k=\p_k^* \SE \ox
L^{\ox k}= \underline{\C} \subset \C^{2n+2} \ox L^{\ox k}$, which
are topologically trivial. We look for everywhere non-zero sections
$s_k$ of $E_k\subset \C^{2n+2} \ox L^{\ox k}$ as they satisfy 
$\p_k=\PP(s_k)$.

Let $\CrP(\f,x)$ be the $C^1$-open property for
sequences of sections $s_k$ of $E_k$ of being $\f$-transverse to
$0$ at the point $x$ (see Definition \ref{skd-trans}). 
We shall use Proposition \ref{local_trans} to 
find sequences of sections $s_k$ which are $\eta$-transverse to
$0$, for some $\eta>0$. Fix any asymptotically holomorphic sequence
$s_k$ of $E_k$ (e.g.\ the zero sections) which will act as the 
starting point of our perturbation process.
Let $x\in M$. Consider the sections $\sref$ of $L^{\ox k}$ given
by Lemma \ref{localized} and define also the local sections of the
line bundle $\p_k^*\SE \subset \underline{\C}^{2n+2}$,
$$
  \s_k: B_{g_k}(x,c) \to \C^{2n+2},
$$
by setting $\s_k(x)$ any vector of norm $1$ in the direction
defined by $\p_k(x)$ and satisfying
the condition $\nabla_r \s_k(y) \perp \s_k(y)$, for any
$y\in B_{g_k}(x,c)$, where $r$ is the radial vector field from $x$.
This determines $\s_k$ uniquely. The following estimates hold:
 \begin{eqnarray}
  & & |\s_k(y)|=1, ~~~~ |\nabla \s_k(y)|=O(1+d_k(x,y)), \nonumber \\
  & & |\bbd \s_k(y)|=O(k^{-1/2}(1+d_k(x,y))), ~~~~ 
  |\nabla \bbd \s_k(y)|=O(k^{-1/2}(1+d_k(x,y))). \label{parallel_b}
 \end{eqnarray}
The first one follows from $\nabla_r \la\s_k,\s_k \ra=
\la\nabla_r \s_k,\s_k \ra+ \la\s_k,\nabla_r \s_k \ra=0$.
For the second one, write $\nabla \s_k=\nabla \p_k+
\la\nabla \s_k,\s_k\ra  \s_k$, where we identify 
$T_{\p_k(y)} \CP^{2n+1}= [\s_k(y)]^{\perp} \subset \C^{2n+2}$, 
isometrically. We already know that $|\nabla \phi_k|=O(1)$. So
\begin{eqnarray*}
 \nabla_r \la \nabla\s_k,\s_k \ra &= &
 \la\nabla \nabla_r \s_k,\s_k \ra+\la \nabla\s_k,\nabla_r\s_k \ra= \\
 &=& \la\nabla \nabla_r \s_k,\s_k \ra+\la \nabla\s_k,\nabla_r\s_k \ra= \\
 &=&-\la\nabla_r\s_k,\nabla \s_k \ra+\la \nabla\s_k,\nabla_r\s_k \ra = \\
 &=&-\la\nabla_r\p_k,\nabla \p_k \ra+\la \nabla\p_k,\nabla_r\p_k \ra =O(1),
\end{eqnarray*}
The first equality uses that $\nabla$ is the standard derivative for
functions with values in $\C^{2n+2}$, and hence the second derivatives
commute. The second equality follows from $\la\nabla_r \s_k,\s_k \ra=0$.
So we have that $\la \nabla\s_k,\s_k \ra=O(d_k(x,y))$ and hence 
$|\nabla \s_k|=O(1+d_k(x,y))$. The other two cases are worked out analogously.

Now define the application
$$
  f_k=\frac{s_k}{\sref \s_k}: B_{g_k}(x,c) \to \C,
$$
which is asymptotically holomorphic by construction. Using a complex
Darboux chart we trivialize $B_{g_k}(x,c)$ to obtain (scaling the 
coordinates by appropiate universal constants) an application
$\hat f_k:B^+ \to \C$ to which we apply Proposition \ref{local_trans}
to obtain $w_k \in B(0,\d)$ such that $\hat f_k -w_k$ is $\eta$-transverse
to $0$ in $B$, where $\eta=\d (\log (\d^{-1}))^{-p}$. Rescaling and 
passing to the manifold, we have that $f_k - Cw_k$ is $C'\eta$-transverse
to $0$, for some universal constants $C$ and $C'$. Define the sequence of
sections $\tau_{k,x}= -w_k \sref \s_k$ of $E_k$, which is asymptotically 
holomorphic and has Gaussian decay by \eqref{parallel_b}, to get 
a perturbation satisfying the conditions in Proposition \ref{globalizate}.

Thus there exists an asymptotically holomorphic sequence $s_k$ of sections
of $E_k$ which is $\eta$-transverse to $0$, for some $\eta>0$. For $k$
large enough, the zeroes of $s_k$ is a symplectic submanifold representing
the trivial homology class, hence the empty set. So $s_k$ is 
nowhere vanishing and hence $\p_k=\PP(s_k)$.

We have that $s_k$ is an asymptotically holomorphic sequence of sections
of $\C^{2n+2}\ox L^{\ox k}$. Let us check that $s_k$ is $\eta$-projectizable, 
i.e.\ that $|s_k|\geq \eta$ everywhere. Suppose that this is not the case and
take the point $x\in M$ where $|s_k|$ attains its minimum. As $|s_k(x)|<
\eta$, $\eta$-transversality implies that $|\nabla s_k(x)|\geq \eta$. Also
$s_k$ is asymptotically holomorphic, so for $k$ large $\nabla s_k(x):
T_xM \to (E_k)_x$ is surjective. Take $v\in T_x M$ such that 
$\nabla_v s_k(x)= s_k(x)$. Evaluating the equality
$$ 
 \nabla |s_k|^2 = \la \nabla s_k, s_k \ra + \la s_k, \nabla s_k \ra.
$$
at the point $x$ and along the direction of $v$, we obtain
$|s_k(x)|^2=0$, which is impossible since we have already proved that
$s_k$ is nowhere vanishing.

Finally the extension to the one-parameter case is trivial.
\hfill $\Box$

\section{Estimated intersections of symplectic submanifolds} \label{sect_esti}

\subsection{Notions on estimated euclidean geometry} \label{angle_s}
In order to set up the definitions needed in Subsection \ref{projective} 
we state the relevant notions and results on angles between
subspaces of euclidean spaces that we shall need. From now on 
we assume that we are in $\R^n$ equipped with the standard euclidean
inner product, but all the proofs apply to a general finite 
dimensional euclidean space.

The angle between two non-zero vectors $v,w\in \R^n$ is defined as
$$ 
\angle(v,w)= \arccos \left(\frac{\langle v,w \rangle}{|v||w|}\right) 
 \in [0,\pi].
$$
The angle is symmetric and satisfies the classical triangular inequality,
$$ 
\angle(u,w) \leq \angle (u,v) + \angle (v,w), 
$$
for non-zero vectors $u,v, w\in \R^n$. Also 
the angle of a vector $u\not =0$ respect to a subspace $V\not
= \{ 0 \}$ is defined as
$$ 
 \angle(u,V)= \min_{v \in V- \{0 \}} \{ \angle (u,v) \}=\angle (u,v(u)) 
  \in [0,\frac{\pi}{2}], 
$$
where $v:\R^n \to \R^n$ is the orthogonal projection onto $V$, well 
understood that when $v(u)=0$ the angle is $\pi/2$.

\begin{definition}
 The maximum angle of a subspace $U \not = \{ 0 \}$ with respect to a 
 subspace $V\not= \{ 0 \}$ is defined as
 $$ 
 \angle_M(U,V)= \max_{u \in U- \{ 0 \}} \angle (u,V). 
 $$
\end{definition}

Notice that this angle is not in general symmetric. But in the 
case $\dim U= \dim V$ symmetry does hold. This is easily checked by 
constructing an orthogonal transformation permuting the two subspaces. 
Indeed the maximum angle $\angle_M(U,V)$ gives a notion of 
proximity between $U$ and $V$ whenever $\dim U\leq \dim V$.

\begin{lemma} \label{sub_add}
 Given $U,V,W$ non zero-subspaces in $\R^n$ then:
 $$ \angle_M(U,W) \leq \angle_M(U,V)+\angle_M(V,W). $$
\end{lemma}

{\bf Proof.}
 We will denote by $v(u)$ the orthogonal projection of the vector $u$ 
 onto the subspace $V$. In the following inequalities, if $v(u)=0$, 
 we suppose that the angle in which this expression appears is $\pi/2$. 
 We have
\begin{eqnarray*}
  \angle_M(U,W) & = & \max_{ u \in U-\{ 0 \}} \{ \min_{w \in W-\{ 0 \}} 
         \{ \angle(u,w) \} \} \leq \\
   & \leq & \max_{ u \in U-\{ 0 \}} \{ \min_{w \in W-\{ 0 \}} \{ 
         \angle(u,v(u))+ \angle(v(u),w) \} \} = \\
   & = &  \max_{ u \in U-\{ 0 \}} \{ \angle(u,v(u)) + 
          \min_{w \in W-\{ 0 \}} \{ \angle(v(u),w) \} \} \leq  \\
   & \leq & \max_{ u \in U-\{ 0 \}} \{ \angle(u,v(u)) \} + \max_{u\in 
    U-\{ 0 \}} \{ \min_{w \in W-\{ 0 \}} \{ \angle(v(u),w) \} \} \leq \\
   & \leq & \angle_M(U,V) +\max_{v\in V-\{ 0 \}} \{ 
    \min_{w \in W-\{ 0 \}} \{ \angle(v,w) \} \} \leq \\ 
   & \leq & \angle_M(U,V) + \angle_M(V,W).
\end{eqnarray*}
\hfill $\Box$

\begin{definition}
The minimum angle between two non-zero subspaces $U,V$ of $\R^n$ 
is defined as follows:
\begin{itemize}
 \item If $\dim U+\dim V< n$ then $\angle_m(U,V)=0$.
 \item If their intersection is not transversal then $\angle_m(U,V)=0$.
 \item If their intersection is transversal then let $W$ be their 
   intersection. Define $U_c$ as the orthogonal subspace in $U$ to 
   $W$, and $V_c$ in the same way. Then 
   $\angle_m(U,V)= \min_{u \in 
   U_c-\{ 0 \}} \{ \angle(u,V_c) \} \in [0,\pi/2]$.
\end{itemize}
\end{definition}

The definition is symmetric because (in the transversal case)
$$ 
  \angle_m(U,V)= \min_{u \in U_c- \{ 0 \}} \{ \min_{v \in V_c- \{0 \}} 
  \{ \angle(u,v)\} \}
$$
 and the two minima commute. Also $\angle_m(U,V)= \min_{u \in 
 U_c-\{ 0 \}} \{ \angle(u,V) \}$.

\begin{lemma}\label{angle_perp}
 For non-zero subspaces $U$ and $V$ of $\R^n$ we have that
$$
  \angle_m(U,V)= 
  \min_{u \in U^{\perp}-\{ 0 \}} \{ \min_{v \in V^{\perp}-\{ 0\} } 
  \{ \angle(u,v) \} \} 
$$
\end{lemma}

{\bf Proof.}
 This is trivial in the case $\dim U+\dim V<n$ or when $U$ and $V$ 
 do not intersect transversely. In the transversal case, 
 we can restrict ourselves to the subspace $(U\bigcap V)^{\perp}$ 
 to compute the angles. So without loss of generality
 we can suppose that $U\oplus V=\R^n$, $U_c=U$ and $V_c=V$.
 As $\dim U=\dim V^{\perp}$, we may construct an orthogonal 
 transformation $\phi$ permuting $U$ and $V^{\perp}$, i.e.\
 $\phi(U)=V^{\perp}$ and $\phi(V^{\perp})=U$. Therefore also
 $\phi(V)=U^{\perp}$. So
$$
  \angle_m(U,V)= \angle_m (\phi(U),\phi(V)) =\angle_m (V^{\perp},
  U^{\perp})=\min_{u \in U^{\perp}-\{ 0 \}} \{ \min_{v \in V^{\perp}-\{ 0\} } 
  \{ \angle(u,v) \} \},
$$
 which proves the lemma.
\hfill $\Box$

\begin{proposition} \label{vari_min}
  For non-zero subspaces $U, V, W$ of $\R^n$ we have that 
  $$
  \angle_m(U,V) \leq \angle_M(U,W) +\angle_m(W,V).
  $$
\end{proposition}
  
{\bf Proof.}
  By Lemma \ref{angle_perp} we have that
\begin{eqnarray*}
  \angle_m(U,V) &=& \min_{u \in U^{\perp}-\{ 0 \}} \{ 
     \min_{v \in V^{\perp}-\{ 0\} } \{ \angle(u,v) \} \}  \leq \\
  &\leq & \min_{u \in U^{\perp}-\{ 0 \}} \{\angle(u,w) \}+ 
     \min_{v \in V^{\perp}-\{ 0\} } \{ \angle(w,v) \},
\end{eqnarray*}
  for any $w\in \R^n$. Choose $w_0\in W^{\perp}-\{0\}$ satisfying 
$$
  \angle_m (W,V)= \min_{w \in W^{\perp}-\{ 0 \}} 
  \{ \min_{v \in V^{\perp}-\{ 0\} } \{ \angle(w,v) \} \}  
  = \min_{v \in V^{\perp}-\{ 0\} } \{ \angle(w_0,v) \}.
$$
 Then we have
$$
  \angle_m(U,V) \leq 
  \min_{u \in U^{\perp}-\{ 0 \}} \{\angle(u,w_0) \}+\angle_m(W,V) \leq 
  \angle_M (W^{\perp}, U^{\perp}) +\angle_m(W,V).
$$

 The result follows once we show that
 $\angle_M (W^{\perp}, U^{\perp}) =\angle_M (U,W)$.
 Put $\angle_M(U,W)=\a$. Let $u\in U$ with $\angle (u,W)=\a$. 
 Denoting by $w$ the projection of $u$ onto $W^{\perp}$, we have that
 $\angle(u,W^{\perp})=\angle (u,w)=\frac{\pi}2 -\a$. So $\angle (w,U)
 \leq \frac{\pi}2-\a$ and hence $\angle(w,U^{\perp})\geq \a$.
 This implies that $\angle_M (W^{\perp}, U^{\perp}) \geq \a= 
 \angle_M (U,W)$. The opposite inequality follows by symmetry.
\hfill $\Box$

\begin{corollary} \label{vari_min-cor}
  Given non-zero subspaces $U, U', V$ of $\R^n$ with 
  $\angle_m(U,V)>\f$ and $\angle_M(U, U')<\delta$ then $\angle_m(U',V)>
  \f-C\d$, where $C$ is a universal constant ($C=1$ in fact). 
 \hfill $\Box$
\end{corollary}

The following result will be very important for our purposes.

\begin{proposition} \label{geometric}
 Given $\f>0$ and $U\in \Gr(m,n)$, $V\in \Gr(r,n)$ subspaces verifying that 
 $\angle_m(U,V)>\f$. Then there are $\g_0>0$ and a constant $C$, 
 depending only on $\f$, such that for any $\g < \g_0$, 
 if $U'\in Gr(m,n)$ and $V'\in Gr(r,n)$ verify that
$$ 
  \angle_M(U,U') <\gamma, ~~~~ \angle_M(V,V') <\gamma, 
$$
 then $U'$ and $V'$ intersect transversally and 
 $\angle_M(U\bigcap V, U' \bigcap V')< C\g$.
\end{proposition}

{\bf Proof.} 
 By Proposition \ref{vari_min} choosing $\gamma_0>0$ small enough, 
 only depending on $\epsilon$, we can assure that the following
 intersections are transversal $U\bigcap V=W$, $U \bigcap V'$, 
 $U' \bigcap V$ and $U' \bigcap V'=W'$
 and that $\angle_m(U',V')\geq \f/2$. By Lemma \ref{sub_add} we have
$$ 
 \angle_M(W,W')\leq \angle_M(W, U\bigcap V')+\angle_M(W', U\bigcap V'). 
$$
 We are going to bound the first term in the right hand side of the 
 inequality, the bounding of the second term being analogous. 

 Put $s= \dim W=r+m-n$. Choose an orthonormal basis
 $(e_1, \ldots, e_s)$ of $W$, extend it to an orthonormal basis 
 $(e_1, \ldots, e_r)$ of $V$ and finally extend it to an orthonormal 
 basis $(e_1, \ldots, e_n)$ of $\R^n$. Note that $(e_{s+1}, \ldots, e_r)$ 
 is an orthonormal basis of $V_c$. 
 As $\angle_m(U,V)=\angle_m(U_c,V)>\f$ and $\angle_M (V,V')<\g_0$ we have
 $\angle_m(U_c,V')>\f/2$ (decreasing $\g_0$ if necessary). 
 So 
 $U_c \cap V'=\{0\}$. Recalling that $V \oplus U_c=\R^n$, we see that
 there is a basis $(e_1+\e_1, \ldots, e_r+\e_r)$ for $V'$ where 
 $\e_j\in U_c$.
 Using that $\angle_m(U,V)> \epsilon$ and that the decomposition
 $\R^n= W\oplus V_c \oplus V^{\perp}$ is orthogonal, we have
\begin{eqnarray*}
 \pr_W^{\perp}(\e_j) & = & 0, \\
 \pr_{V_c}^{\perp}(\e_j) & \leq & |\cos \epsilon| |\e_j|, \\
 \pr_{V^{\perp}}^{\perp}(\e_j) & \geq & \sqrt{1-|\cos \epsilon|^2}
  |\e_j|=|\sin\f| |\e_j|. 
\end{eqnarray*}
Checking the angle of $e_j+\e_j$ with respect to $V$, we get that
\begin{equation}
  \angle_M(V,V') \geq \arctan \frac{|\sin \epsilon| |\e_j|}{1+|\cos \f|
  |\e_j|} \geq \arctan \left( \frac{\sin \epsilon}{1+|\e_j|} |\e_j| \right). 
  \label{epsi_cotas}
\end{equation}
 For $\g_0< \arctan \frac{\sin\f}{2}$, \eqref{epsi_cotas} implies that 
 $|\e_j|<1$ and hence
 we get that $\angle_M(V,V')\geq \arctan (\frac{\sin\f}{2} |\e_j|) \geq
 \frac{4}{\pi}\frac{\sin\f}{2}|\e_j|$, or said otherwise $|\e_j| <C \angle_M
 (V,V')$ for a constant $C$ depending on $\f$.

 Now let us compute $\angle_M(W,U\bigcap V')$. The intersection $U\bigcap
 V'$ has basis $(e_1+\e_1, \ldots, e_s+\e_s)$. Take a general vector
 $u=\sum_{i=1}^s a_i(e_i+\e_i)$ in $U\bigcap V'$ and compute $\angle (u,W)$.
 We may suppose that $a=(a_1,\ldots, a_s)$ has norm one. 
 Write $\e=\sum_{i=1}^s a_i \e_i$. Then
$$
 \angle (u,W)= \arccos \frac{1}{\sqrt{1+|\e|^2}}=\arctan |\e| \leq 
 |\e|.
$$
 Finally 
$$
 \angle_M(W ,U\bigcap V') \leq \max_{|a|=1} |\sum_{i=1}^s a_i \e_i|= 
 \max_{1\leq i\leq s} |\e_i| \leq C \angle_M(V,V') \leq C \g.
$$
\hfill $\Box$

Now we are going to set up the relationship between the transversality
of maps in the Donaldson-Auroux approach and the angles defined 
above. This is the content of the following

\begin{lemma}\label{bridge}
 Let $U, V$ be two non-zero subspaces of $\R^n$ and let 
 $g:U \to V$ and $h:U \to V^{\perp}$ be the projections from $U$ with
 respect to the decomposition $\R^n=V\oplus V^{\perp}$.
 If $h$ has a right inverse $\theta$ satisfying 
 $|\theta|<\g^{-1}$ for some $\g>0$ then $\angle_m(U,V)> \g$. 
\end{lemma}

{\bf Proof.}
 In the first place, as $h$ is onto, the intersection between $U$ and
 $V$ is transversal. Let $W=U\cap V$. Define $\hat{\theta}=
 \pr^{\perp}_{U_c} \circ \theta:V^{\perp} \to U_c$, 
 which is an inverse of $h: U_c\to V^{\perp}$
 such that $|\hat{\theta}|<\g^{-1}$. Now consider any 
 $u \in U_c-\{0\}$ and put $v=h(u)$. Then 
$$
 \angle(u,V)=\arcsin \frac{|h(u)|}{|u|} = 
 \arcsin\frac{|v|}{|\hat{\theta}(v)|} > 
 \arcsin \frac{1}{\gamma^{-1}} > \g, 
$$
 and the proof is concluded. 
\hfil $\Box$

\subsection{Projective symplectic geometry} \label{projective}
In this Subsection we will prove Theorem \ref{good_inter0}. This will
provide a geometric proof of Bertini's theorem, the 
main result of \cite{Do96}. Although our proof is more technical 
and long, it has the advantage of giving us a more general kind 
of symplectic submanifolds than those in \cite{Do96,Au97}. In fact our
technique will allow us a simple generalization to solve the problem 
of constructing determinantal symplectic submanifolds in 
Section \ref{determinantal}. First of all, in order to measure 
the holomorphicity of submanifolds, let us introduce the complex angle 
of even dimensional subspaces $V\subset \C^n$ as
\begin{eqnarray*}
  \beta: \Gr_{\R}(2r,2n) & \to & [0,\pi/2] \\
  V & \to & \angle_M(V,J V).
\end{eqnarray*}
 Clearly $\b(V)=0$ if and only if $V$ is complex and 
 $\beta(V)<\pi/2$ if and only if $V$ is symplectic.

\begin{definition} \label{def_ah}
 Let $(M,\o)$ be a symplectic submanifold endowed with a compatible
 almost complex structure $J$. 
 A sequence of submanifolds $S_k \subset M$ is asymptotically holomorphic
 if $\b(T S_k)=O(k^{-1/2})$.
\end{definition}

Note that if $S_k$ are asymptotically holomorphic submanifolds 
then they are symplectic for $k$ large.
If $\p_k:M\to \CP^N$ is a sequence of asymptotically holomorphic 
embeddings then $\p_k(M)$ is a sequence of 
asymptotically holomorphic submanifolds.

\begin{proposition} \label{angle}
 Let $\phi_k^1:(M_1, J_1) \to \CP^N$ and $\phi_k^2:(M_2, J_2) \to \CP^N$
 be two sequences of asymptotically holomorphic embeddings. 
 Suppose that there exists $\epsilon>0$ independent of $k$ such that 
 for any  $x \in \phi_k^1(M_1) \bigcap \phi_k^2(M_2)$, the minimum angle
 between $(\phi_k^1)_* TM_1 (x)$ and $(\phi_k^2)_* TM_2 (x)$ is greater
 than $\epsilon$. Then $S_k=\phi_k^1(M_1) \bigcap  \phi_k^2(M_2)$ is a
 sequence of 
 asymptotically holomorphic submanifolds (hence symplectic for $k$ large).
 Also $S^j_k=(\phi_k^j)^{-1}(S_k)$ is a sequence of asymptotically
 holomorphic submanifolds of $M_j$, $j=1,2$. Moreover there exists a sequence
 of compatible almost complex structures $J_k^j$ of $M_j$ such that $S_k^j$ is
 pseudoholomorphic for $J_k^j$, $|J_k^j-J_j|=O(k^{-1/2})$ and $\phi^j_k$ 
 restricted to $(S^j_k, J_k^j)$ is a sequence of asymptotically 
 holomorphic embeddings in $\CP^N$, $j=1,2$.

 The same statement holds for the case of one-parameter families of
 embeddings $(\phi_{t,k}^1)_{t \in [0,1]}$ and $(\phi_{t,k}^2)_{t \in [0,1]}$.
\end{proposition}

 Remark that $M_1$ and $M_2$ are not necessarily compact manifolds. 

\noindent {\bf Proof.} 
 Let $J_0$ be the standard complex structure of $\CP^{2n+1}$. Then
 $\angle_M((\phi_k^j)_* TM, J_0 (\phi_k^j)_* TM)=O(k^{-1/2})$ for $j=1,2$.
 By Proposition \ref{geometric}, $\angle_M(TS_k, J_0 TS_k)=O(k^{-1/2})$. 
 As $|(\p_k^j)_*J_j-J_0|=O(k^{-1/2})$ on $(\p_k^j)_*TM$, we have
 $\angle_M(TS_k, (\p_k^j)_* J_j TS_k)=O(k^{-1/2})$ and so
 $\angle_M(TS_k^j, J_j TS_k^j)=O(k^{-1/2})$, i.e.\ $S_k^j$ is a sequence of
 asymptotically holomorphic submanifolds of $M_j$.

 Finally we have to build $J_k^j$ on $M_j$ such that $|J_k^j-J_j|=O(k^{-1/2})$ 
 and $S_k^j$ is $J_k^j$-holomorphic. Take the composition $\tilde{J}_k^j:
 TS_k^j \subset TM 
 \stackrel{J_j}{\to} TM \stackrel{\pr^{\perp}}{\surj} TS_k^j$ with
 square close to $-1$, for $k$ large enough. So we can homotop it to 
 an almost complex structure $J_k^j$ on $S_k^j$. Then we extend this
 $J_k^j$ to a small tubular neighborhood of $S_k^j$ by giving a complex
 structure to the normal bundle of $S_k^j$. Finally a homotopy between
 $J_k^j$ and $J_j$ allows us to extend $J_k^j$ off a little bigger
 neighborhood of $S_k^j$ matching with $J_j$ on the border. This gives
 the required $J^j_k$.

 The result for continuous one-parameter families is trivial from the 
 non-parametric case.
\hfill $\Box$

Let us have a smooth submanifold $N$ of a manifold $X$.  
If we fix a metric on $X$ we can define a geodesic flow $\varphi_t$.
In particular, following the perperdicular directions to $N$ we can 
identify a tubular neighborhood of the zero section of the normal
bundle of $N$ (defined as $|n|<t_0$, $n \in \nu(N)$, for some 
small $t_0>0$) with a tubular neighborhood $U_N\subset X$ of $N$. 
So we can define an integrable distribution $D_N$ in $U_N$ as
$$ 
 D_N(\varphi_n(x))=(\varphi_n)_* T_xN, \ \forall x\in N, n\in \nu(N), 
 |n|<t_0.
$$
where $(\varphi_n)_*$ denotes parallel transport along the geodesic 
tangent to $n$. 

\begin{definition} \label{openangle}
 Suppose $\phi_k:M\to X$ is a sequence of asymptotically holomorphic 
 embeddings into a Hodge manifold $X$.
 Let us fix a complex submanifold $N\subset X$. We say that $\p_k$ is
 $\s$-transverse to $N$, with $\s<t_0$, if for all $x\in M$ and all $k$, 
$$ 
  d(\p_k(x),N)< \sigma \Rightarrow \angle_m((\phi_k)_*(T_xM),D_N(\p_k(x))) 
  >\s.
$$
\end{definition}

 This property is $C^1$-open, i.e.\
 given $\phi_k$ an embedding $\eta$-transverse to $N$, then a perturbation of
 $\hat\p_k$ with $d_{C^1}(\p_k,\hat\p_k) <\d$ 
 is $(\eta-C\d)$-transverse to $N$, where $C$ is a universal constant.

 Obviously a $\s$-transverse sequence of embeddings $\p_k$
 verifies the conditions of Proposition \ref{angle} with 
 $\p_k^1=\p_k:M \to X$ and $\p_k^2=i:N\inj X$. 
 The following result then completes the proof of Theorem \ref{good_inter0}

\begin{theorem} \label{good_inter1}
 Let $\phi_k= \PP (s_k)$, where $s_k$ is an asymptotically holomorphic
 sequence of sections of $\C^{2n+2} \ox L^{\ox k}$ which
 is $\g$-projectizable and $\g$-generic of order $n$, for some $\g>0$.
 Let us fix a holomorphic submanifold $N$ in $\CP^{2n+1}$. Then 
 for any $\d>0$ there exists an asymptotically holomorphic sequence of
 sections $\s_k$ of $\C^{2n+2} \ox  L^{\ox k}$ such that
\begin{enumerate}
  \item $|\s_k-s_k|_{g_k,C^1}<\delta$.
  \item $\hat\p_k=\PP (\s_k)$ is a $\eta$-asymptotically holomorphic embedding
    in $\CP^{2n+1}$ which is $\f$-transverse to $N$, for some $\eta>0$ and
    $\epsilon>0$. In the case $\dim M + \dim N < 2n+1$ we actually have 
    that $d_{FS}(\hat\p_k(M),N))>\f$, for $k$ large enough.
\end{enumerate}

 Moreover the result can be extended to one-parameter continuous families 
 of complex submanifolds $(N_t)_{t\in [0,1]}$, taking in this case as 
 starting point a continuous family $\phi_{t,k}=\PP(s_{t,k})$ where $s_{t,k}$
 are asymptotically $J_t$-holomorphic sections of $\C^{2n+2} \ox L^{\ox k}$
 which are $\g$-projectizable and $\g$-generic of order $n$, for some $\g>0$.
\end{theorem}

 The proof of this result will be the content of Subsection \ref{esti_inter}.
 Now we shall extract some corollaries from it. 
 The first one is the main Theorem of \cite{Do96}.

\begin{corollary} \label{bertini}
 Given a compact symplectic manifold $(M,\o)$, suppose that 
 $[\o/2\pi]\in H^2(M, \R)$ is the reduction of an integral class $h$. 
 Then for $k$ large enough there exists symplectic submanifolds realizing 
 the Poincar\'e dual of $kh$. Moreover, perhaps by increasing $k$,
 we can assure that all the symplectic submanifolds realizing this 
 Poincar\'e dual,
 constructed as transverse intersections with a fixed complex hyperplane 
 of asymptotically holomorphic sequences of embeddings with
 respect to two compatible almost complex structures, are isotopic.
 The isotopy can be made by symplectomorphisms.
\end{corollary}

Recall that we obtain an isotopy result similar to \cite{Au97}, where the 
isotopy of the submanifolds obtained as zero sets of a special
set of sections of the line bundle $L^{\otimes k}$ is obtained.
The Auroux' more general case of vector bundles will be proved in 
Section \ref{emb_grass}.

{\bf Proof.}
 The existence result is a direct consequence of the previous statements. By 
 Theorem \ref{existence1} we build an asymptotically holomorphic sequence of
 embeddings to  $\CP^{2n+1}$. In $\CP^{2n+1}$ we choose a complex 
 hyperplane $H$. By Theorem \ref{good_inter1}
 we perturb the sequence of embeddings to find a new asymptotically
 holomorphic sequence of embeddings $\phi_k$ such that $\phi_k(M)$ intersects
 $H$ with minimum angle greater than $\epsilon>0$. Finally 
 using Proposition \ref{angle} we
 obtain that $\phi_k(M) \bigcap H=H_M$ is an asymptotically holomorphic
 sequence of submanifolds, and these manifolds are symplectic
 for $k$ large enough. Also $\phi_k^{-1}(H_M)$ is a symplectic
 submanifold of $M$ for $k$ large enough. A direct topological argument 
 shows us that it is Poincar\'e dual of $kh$.

 For the isotopy statement, let us assume that there are two 
 sequences of symplectic 
 submanifolds $W_k^0$ and $W_k^1$, both Poincar\'e dual of $kh$, obtained 
 as intersections between two $\eta$-asymptotically
 $J_j$-holomorphic sequences $\PP(s_{k,j})$, $j=0,1$, 
 and two fixed complex hyperplanes $H_0$ and $H_1$ in $\CP^{2n+1}$ with
 angles greater than a fixed $\epsilon>0$. Then we will prove that in this 
 case they are isotopic. We only have to construct the straight segment 
 $H_t$, in the dual space, of hyperplanes connecting
 $H_0$ and $H_1$. Also we define the following family of asymptotically 
 holomorphic sequences:
$$ 
  s_{t,k}=\left\{ \begin{array}{lll} (1-3t)s_{0,k}, \qquad & \text{with } 
          J_t=J_0, & t\in [0,1/3] \\
        0, & \text{with } J_t=\text{Path}(J_0,J_1), \qquad & t\in [1/3,2/3] \\ 
        (3t-2)s_{1,k}, & \text{with } J_t=J_1, & t\in [2/3,1].
  \end{array} \right. 
$$
 By means of Theorem \ref{good_inter1}, we obtain a family $\phi_{t,k}
 =\PP(\s_{t,k})$ of asymptotically $J_t$-holomorphic embeddings which are
 $\eta/2$-transverse to $N$, choosing the perturbation $\delta>0$ in the
 statement of the theorem, in such a way that
\begin{equation}\label{not_destroy} 
  \eta-C\delta >\eta/2, 
\end{equation} 
 where $C$ is the universal constant of the $C^1$-openness 
 of the transversality to $N$. This gives us a family of symplectic 
 isotopic submanifolds $(W_k^t)'$ in $M$ for each fixed large $k$. 
 The problem is that $W_k^0$ does not coincide with $(W_k^0)'$ (and 
 respectively for $t=1$). Using \eqref{not_destroy} we can assure that
 they are isotopic, in fact the linear segment 
 $((1-t)\sigma_{0,k}+ts_{0,k})_{t\in[0,1]}$ provides a family 
 of asymptotically holomorphic embeddings transverse to $H_0$, for $k$ large
 enough giving the desired isotopy. 
\hfill $\Box$

The constructive technique of Theorem \ref{good_inter1} is more general 
because we do not have to choose hyperplanes in 
$\CP^{2n+1}$ to make the intersection. However, the difficulty in finding
topological information about the constructed submanifolds makes that we
cannot assure that they are more general that the ones produced
in \cite{Au97}. To overcome this problem we are going to construct in
Section \ref{determinantal} a special
kind of submanifolds where we can compute symplectic invariants
using similar results from algebraic geometry.

\subsection{Estimated intersections in $\CP^{2n+1}$.} \label{esti_inter}
 Now we aim to prove Theorem \ref{good_inter1}. 
 Our objective is to find sequences $\p_k$ of asymptotically holomorphic
 embeddings which are $\s$-transverse to $N$.

\noindent {\bf Proof of Theorem \ref{good_inter1}.}
As usual we begin with the simplest case, when the complex codimension
of $N$ is $1$. Also we consider the non-parametric case, being the parametric 
one a simple generalization. We say that a sequence of sections $s_k$ which is
$\g/2$-projectizable and $\g/2$-generic of order $n$ verifies $\CrP(\f,x)$ 
if $\PP(s_k)$ is $\f$-transverse to $N$ at the point $x$. This property 
is local and open in $C^1$-sense, for $\f<t_0$. To make use of 
Proposition \ref{globalizate} we need to find local sections with Gaussian 
decay obtaining local transversality.
To achieve this local transversality we are going to use Proposition
\ref{local_trans}. (We could have used instead the case $m=1$ proved 
in \cite{Do96,Au97}, by increasing a little the complications of the 
globalization process, which is the way followed by Auroux 
in \cite{Au97,Au99}.)

As $N$ is a fixed holomorphic submanifold, we may fix a finite covering of
$\CP^{2n+1}$ by balls $U_j$ such that $N$ is defined as the zero set of a 
holomorphic function $f_j:U_j \to \C$ in each $U_j$ and such that for
any $z_1, z_2 \in U_j \cap U_N$, $\angle_M (D_N(z_1),D_N(z_2)) \leq \e$,
and for any $z_1, z_2 \in U_j$, $\angle_M (\ker df_j(z_1),\ker df_j(z_2)) \leq \e$,
with $\e>0$ an arbitrarily small number fixed along the proof.

We choose a constant $C$ independent of $k$ such that 
$|\nabla \p_k|_{g_k} \leq C$. Therefore
$\p_k(B_{g_k}(x,c)) \subset B_{g_{FS}}(\p_k(x),Cc)$, for any $c$.
Now we choose $c>0$ small enough satisfying the following premises:
\begin{enumerate}
\item 
 Let $x\in M$. With a transformation of $U(2n+2)$ in $\C^{2n+2}$, we may 
 suppose that $s_k(x)=(s_k^0(x), 0,\ldots, 0)$. As $s_k$ is 
 $\g$-projectizable and asymptotically holomorphic, we can choose a 
 universal $g_k$-radius $c$ with $|s_k^0| \geq \g/2$ on $B_{g_k}(x,20c)$.
 Also the sections $\sref$ of Lemma \ref{localized} satisfy $|\sref|\geq
 c_s$ on $B_{g_k}(x,20c)$. Note that $\p_k(B_{g_k}(x,20c)) 
 \subset B_{g_{FS}}(\p_k(x),20Cc)$.

\item 
 We use the standard chart $\P_0$ for $\CP^{2n+1}$ around 
 $p=\p_k(x)=[1,0,\ldots ,0]$ to trivialize the ball 
 $B_{g_{FS}}(p,20Cc)$. We may choose $c$ small enough so that
 $\P_0$ is near an isometry, in the sense that
 $$
  \frac23 |\P_0(q)| \leq d_{FS}(p,q) \leq 2 |\P_0(q)|.
 $$
 for $q \in B_{g_{FS}}(p,20Cc)$. Also we require $|\nabla \P_0|\leq 2$
 in such ball. With respect to this trivialization 
 the map $\p_k$ is given locally as
 \begin{eqnarray*}
  f_k=\P_0\circ \p_k : B_{g_k}(x,20c) &\to & B(0, 40Cc) \\
  y &\mapsto &\left( \frac{s_k^1(y)}{s_k^0(y)}, \ldots, 
  \frac{s_k^{2n+1}(y)}{s_k^0(y)} \right).
 \end{eqnarray*}
 Clearly $|\nabla f_k|\leq 2C$ uniformly in $k$.
\item 
 We can reduce $c$ so that, for any $p$, 
 $B_{g_{FS}}(p,20Cc) \subset U_j$ for some $U_j$. 
 Therefore $N$ is defined in 
 $B(0,15Cc)$ by a function $f:B(0,15Cc)\to \C$. Call $Z=Z(f)$ in such
 ball. The angle condition means that $\ker df(z_1)$, $\ker df(z_2)$ are close
 enough (say less than $\pi/6$) for $z_1,z_2\in Z$.
\end{enumerate}

Let $x\in M$. 
In the case $d(\phi_k(x),N) \geq 2Cc$, as we perform a small perturbation, say
of norm $\d>0$ such that $d_{FS}(\phi_k(x),\hat{\phi}_k(x))< \frac12 
Cc$, for all $x\in M$, there is still $\frac12 Cc$-transversality at 
a $c$-neighbourhood of $x$. So we are finished.

Suppose $d(\phi_k(x),N)< 2 Cc$. Then take a point $z_0 \in B(0,4Cc)\cap Z$ 
which gives the minimum distance from $0$ to $Z$. If $0 \notin Z$, take
$v=(v_1,\ldots, v_{2n+1}) \in \C^{2n+1}$ a unitary vector in the 
direction of the complex line from $0$ to $z_0$. 
This vector is perpendicular to $T_{z_0}Z$. If $0\in 
Z$ then let $v$ be a unitary vector orthogonal to $T_0Z$. Therefore
\begin{equation} \label{angle_f}
 \la df(z),v \ra\geq \frac12 |df(z)|
\end{equation}
for any $z\in Z\cap B(0,15 Cc)$, by the condition on the angle
(taking $\f>0$ small enough).

Let $r_0\in \C$ with $r_0 v=z_0 \in Z$. We look for a function 
$r_k=r_k(y):B_{g_k}(x,c) \to \C$ such that $r_k(x)=r_0$ and 
\begin{equation} \label{ift}
  f \left( f_k^1(y)+r_k v_1, \ldots, f_k^{2n+1}(y)+r_k v_{2n+1} \right)=0.
\end{equation}
This corresponds to tracing a straight line from the image of the point $y\in
B_{g_k}(x,c)$ to $Z$ with direction $v$. Such $r_k$ can be found with the 
use of the implicit function theorem applied to the function 
$F:B_{g_k}(x,c) \x B(r_0, 4Cc) \to \C$ given as the left hand side 
of \eqref{ift}.
This $F$ is well-defined since $f$ is defined on $B(0,10Cc) \subset \P_0(U_j)$.
To guarantee the existence of $r_k=r_k(y)$ for all $y\in B_{g_k}(x,c)$
we have to check that
$$
 \left|\frac{\nabla_y F}{\bd F/\bd r_k}\right|=
 \left|\frac{\la df,\nabla f_k \ra}{\la df, v\ra}\right|
 \leq 4C,
$$
which holds thanks to \eqref{angle_f}.
This gives the existence of $r_k$ in the whole of the ball $B_{g_k}(x,c)$
as well as the bound $|\nabla r_k|\leq 4C$, and hence $|r_k|\leq 8Cc$.

Now our task will be to prove that $r_k$ is asymptotically holomorphic, so
we change a geometrical transversality problem into a local one.
For this let us compute $\bbd r_k$.  
Recall that $f_k$ is asymptotically holomorphic and $f$ is holomorphic.
Differentiate the equality $f(f_k(y) + r_k (y) v)=0$ to get 
\begin{eqnarray} \label{inter}
 0&=& \bbd (f(f_k(y)+r_k(y)v)) = \bd f (z) \cdot (\bbd f_k(y) +
  \bbd r_k(y) v)= \nonumber \\ &=&  O(k^{-1/2})+ \la df(z),v \ra \bbd r_k(y),
\end{eqnarray}
with $z=f_k(y)+r_k(y)v$.
Using \eqref{angle_f} we get that $\bbd r_k = O(k^{-1/2})$. We
already know that $|\nabla r_k| = O(1)$.
Differentiating \eqref{inter} one easily obtains also that
$|\nabla \bbd r_k|=O(k^{-1/2})$. So $r_k$ is asymptotically holomorphic. 
We shall achieve transversality for the function
$$
  h_k=r_k\frac{s_k^0}{\sref}: B_{g_k}(x,c)\to \C,
$$
which is also asymptotically holomorphic.

Dividing $h_k$ by an appropriate constant,
using the chart $\Phi_k$ defined in Lemma \ref{trivializa} and scaling the
coordinates by a universal constant, we obtain a function $\tilde{h}_k$
defined on $B^+$ satisfying the hypothesis of Proposition \ref{local_trans},
for $k$ large enough. So going back to $h_k$ through universal constants,
we find $|w_k|<\delta$ such 
that $h_k-w_k$ is $\eta$-transverse to $0$ with $\eta=c'\delta(\log 
(\delta^{-1}))^{-p}$.

Now we have a direction $v$ and a modulus $w_k$ for a perturbation. 
The perturbation we give is
$$ 
 \tau_{k,x}= (0, -w_kv_1 \sref, \ldots, -w_k v_{2n+1} \sref).
$$
Let us look at the perturbed map $\hat{\p}_k=\PP(s_k+\tau_{k,x})$. 
It is asymptotically holomorphic and $\g'$-projectizable and $\g'$-generic 
of order $n$, for some $\g'>0$, with $|\tau_{k,x}|<c''\d$
(for $\d>0$ small enough). Let us check that $\hat{\p}_k$ is 
$\eta$-transverse to $N$ with $\eta=
c'\delta (\log (\delta^{-1}))^{-p}$ and $c'$ a constant depending only on 
$c$ and the asymptotically holomorphic bounds of $s_k$. With this, 
applying Proposition \ref{globalizate}, the proof in this case is concluded. 
Only a little problem may appear, that the deformed embedding can become an 
immersion, but then an arbitrarily small perturbation solves the problem.

 The $h_k$ associated to $\hat{\p}_k$ is $\hat{h}_k=h_k -w_k$.
 The final point is to set up the relationship between the transversality of 
 $\hat{h}_k$ to $0$ and the transversality of $\hat{\phi}_k$ to $N$. Note that 
 we have $\hat{r}_k=\hat{h}_k \frac{\sref}{s^0_k}$,
 $\hat{f}_k=\P_0 \circ \hat{\p}_k$ and
 $\hat{\pi}_k=\hat{f}_k+\hat{r}_k v=\pi_k$.

Using that $|\sref/s_k^0|$ is bounded above and below uniformly and that
$|\nabla (\sref/s_k^0)|=O(1)$, it is easy to prove that if $\hat{h}_k$ is
$\eta$-transverse to $0$ then $\hat{r}_k$ is $c_0\eta$-transverse to $0$,
for some universal constant $c_0$.

 Let $y\in B_{g_k}(x,c)$. If $|\hat{r}_k(y)| \geq c_0\eta$ 
 then $d(\hat{\p}_k(y),N)\geq c_1\eta$, for some universal constant $c_1$.
 Otherwise $|\nabla \hat{r}_k(y)| >c_0\eta$.
 We shall use Lemma \ref{bridge}
 for the subspaces $U=(d\hat{f}_k)_*T_y M$ and $V=T_{\pi_k(y)}Z$
 of $\C^{2n+1}$. Let $V'=[v]$. 
 The projections from $U$ to the summands of the decomposition 
 $\C^{2n+1}=V\oplus V'$
 are given respectively by $g=d\pi_k\circ (d\hat{f}_k)^{-1}$ and 
 $h=-v\, d\hat{r}_k\circ (d\hat{f}_k)^{-1}$. This follows from
 $d\pi_k=d\hat{f}_k +d\hat{r}_k\, v $ which gives
 $\text{Id}=d\pi_k\circ (d\hat{f}_k)^{-1}-v\,
 d\hat{r}_k\circ (d\hat{f}_k)^{-1}$.
 The map $h$ has a right inverse of norm
 bounded by $C'\eta^{-1}$, for some universal constant $C'$ (here we use that 
 $\phi_k$ is generic of order $n$ and that the perturbations are small).
 It is easy to check that Lemma \ref{bridge} is still valid when $V$ and 
 $V'$ are almost orthogonal (and not just orthogonal), so we have
$$
  \angle_m((d\hat{f}_k)_*T_y M, T_{\pi_k(y)}Z) \geq c_2 \eta.
$$ 

 Push forward the distribution $D_N$ through
 the chart $\P_0$ to a distribution $D_Z$ in $B(0,15Cc)$.
 Then there exists a constant $C''$ independent of $k$ such that
$$
  \angle_M(T_z Z, D_Z(z+\l v)) <C'' d(z+\l v,Z),
$$
 for $z\in Z$, $\l \in \C$ with $|z| < 14Cc$, $|\l| <Cc$.
 Now use Proposition \ref{vari_min} to get 
$$
  \angle_m((d\hat{f}_k)_*T_y M,D_Z(\hat{f}_k(y)))>c_2\eta-C''d(\hat{f}_k(y),Z).
$$
 For $d(\hat{f}_k(y),Z)< c_2\eta/2C''$ we get 
 $\angle_m ((d\hat{f}_k)_*T_y M,D_Z(\hat{f}_k(y))) > c_2 \eta/2$.
 Passing to the manifold we get
 $\angle_m((d\hat{\p}_k)_*T_y M, D_N (\hat{\p}_k(y)) > c_2'\eta$, 
 whenever $d(\hat{\p}_k(y),N)< c_1'\eta$, 
 for some universal constants $c'_1$ and $c'_2$.

To achieve the solution when the codimension of $N$ is $r>1$, we follow 
the same ideas than in the precedent case. 
In this case $f:B(0,15Cc) \to \C^r$ and one chooses the point $z_0$ 
giving the minimum distance from $0$ to $Z$
which yields a vector $v_1$ orthogonal to $Z$ at $z_0$. Then one completes
to an unitary basis $(v_1, \ldots, v_r)$ for the orthogonal to $T_{z_0}Z$.
The function $r_k: B_{g_k}(x,c) \to \C^r$ is defined by the
condition $f(f_k+r_k^1v_1 +\ldots +r_k^rv_r)=0$. The perturbation will be 
of the form
$$ 
  \tau_{k,x}=-(0, w_k^1 v_1^1\sref+\cdots +w_k^r v_r^1\sref,
  \ldots, w_k^1v_1^{2n+1}\sref+\cdots+w_k^r v_r^{2n+1}\sref), 
$$
where $v_i=(v_i^1, \ldots, v_i^{2n+1})$, $i=1,\ldots, r$ and
$w_k=(w_k^1, \ldots, w_k^r)\in \C^r$. The proof above works out in this case.
\hfill $\Box$

\section{Asymptotically holomorphic embeddings to grassmannians} 
\label{emb_grass}

Let $(M,\omega)$ be a symplectic manifold of integer class and let $L$ 
stand for the hermitian line bundle with a connection $\nabla$ with 
curvature $-i\o$. Let $E$ be a rank $r$ hermitian bundle over $M$ endowed
with an hermitian connection. Fix a compatible almost complex structure
$J$ on $M$. In this Section we shall deal with the issue of constructing 
sequences of embeddings of $M$ into the grassmannian $\Gr(r,N)$ which are 
asymptotically $J$-holomorphic in the sense of Definition \ref{def:ah}.
More specifically, we aim to prove the following result from which
Theorem \ref{existence20} follows.

\begin{theorem} \label{existence21}
 Suppose $N>n+r-1$ and $r(N-r)> 2n$. 
 Given an asymptotically $J$-holomorphic sequence of sections $s_k$ of 
 the vector bundles $\C^{N}\ox E \otimes L^{\otimes k}$ and $\alpha>0$ 
 then there exists another sequence $\sigma_k$ verifying that:
\begin{enumerate}
  \item $|s_k-\sigma_k|_{C^1,g_k}< \alpha$.
  \item $\p_k=\Gr(\s_k)$ is an asymptotically holomorphic sequence of 
    embeddings in $\Gr(r,N)$ for $k$ large enough.
  \item $\phi_k^*\, \SU = E\otimes L^{\otimes k}$, where $\SU \to \Gr(r,N)$ is
    the universal rank $r$ bundle over the grassmannian.
\end{enumerate}
 Moreover given two asymptotically holomorphic sequences $\phi_k^0$ and 
 $\phi_k^1$ of embeddings in 
 $\Gr(r,N)$ with respect to two compatible almost complex structures, then
 for $k$ large enough there exists an isotopy of asymptotically 
 holomorphic embeddings $\phi_k^t$ connecting $\phi_k^0$ and $\phi_k^1$.
\end{theorem}

\subsection{Proof of main result.}
First let us fix some notation.
A point $s\in \Gr(r,N)$ corresponds to an $r$-dimensional subspace 
$V_s \subset \C^N$. Choosing a basis $s_1,\ldots, s_r$ for $V_s$,
we denote 
$$
s=\left[\begin{array}{c} s_1 \\ \vdots \\ s_r \end{array}\right]=
  \left[\begin{array}{ccccc} s_{11} & s_{12} & \cdots & s_{1N} \\
  \vdots & & \ddots & \vdots \\
  s_{r1} & s_{r2} & \cdots & s_{rN} \end{array} \right].
$$
This identifies $s$ as the equivalence class of $r\x N$ matrices of
rank $r$ under the action of $\GL(r,\C)$ on the left.
The standard metric $g_{Gr}$ for $\Gr(r,N)$ is the metric induced by the
Fubini-Study metric $g_{FS}$ under the Pl\"ucker 
embedding \cite[Chapter 1, Section 5]{GH78}
\begin{eqnarray*}
  \Gr(r,N) &\too & \PP(\bigwedge\nolimits^r \C^N) \\
  \left[\begin{array}{c} s_1 \\ \vdots \\ s_r \end{array}\right] 
  & \mapsto & s_1 \wedge \cdots \wedge s_r.
\end{eqnarray*}

We proceed by steps to obtain asymptotically holomorphic embeddings.

\begin{definition}
 Let $\g>0$ and $0\leq l\leq r$. A sequence of asymptotically 
 $J$-holomorphic sections $s_k=(s_k^1,\ldots, s_k^N)$ of 
 the vector bundles $\C^{N} \otimes E \otimes L^{\otimes k}$ is said to be
 $\gamma$-grassmannizable of order $l$ if for all $x\in M$, 
 $|\bigwedge^l s_k(x)|>\gamma$. It is $\g$-grassmannizable when it is
 $\g$-grassmannizable of order $r$. (Here $s_k=(s_k^1,\ldots, s_k^N)$ is
 interpreted as a morphism of bundles
 $\underline{\C}^{N} \to E \otimes L^{\otimes k}$
 and $\bigwedge^l s_k$ is the corresponding $l$-fold wedge product.)
\end{definition}

If we have the condition of $\gamma$-grassmannizability
for a section $s_k$ then we obtain a morphism $\p_k=\Gr(s_k):
M\to \Gr(r,N)$, called the grassmannization of $s_k$, as
follows. At a point $x$ take a basis $(e_1,\ldots, e_r)$ for the
fibre of $E$ at $x$. Then 
$$
 \p_k(x)=\left[ s_k^1(x), \ldots, s_k^N(x) \right]=
  \left[\begin{array}{ccccc} s_k^{11} & s_k^{12} & \cdots & s_k^{1N} \\
  \vdots & & \ddots & \vdots \\
  s_k^{r1} & s_k^{r2} & \cdots & s_k^{rN} \end{array} \right].
$$
where $s_k^i(x) =s_k^{1i}e_1 +\cdots +s_k^{ri}e_r$. This is well-defined
and independent of the chosen basis.

\begin{definition}
 Let $\eta>0$ and $0\leq l\leq n$.
 A sequence of asymptotically $J$-holomorphic $\gamma$-grassmannizable
 sections $s_k$ of vector bundles $\C^N\otimes E\otimes L^{\otimes k}$ 
 is $\eta$-generic of order $l$, with $\eta>0$, if given $\Gr(s_k)$ 
 then for all $x\in M$, $|\bigwedge^l \partial \Gr(s_k)(x)|_{g_k}>\eta$.
\end{definition}

In order to prove Theorem \ref{existence21} we shall 
use the following auxiliar Proposition that 
will be proved in the following Subsections. Also we state the analogue
of Lemma \ref{lem:emb->sec} which will be proved in 
Subsection \ref{lifting2}.

\begin{proposition} \label{key2}
 Suppose $N>n+r-1$ and $r(N-r)> 2n$. 
 Let $s_k$ be an asymptotically $J$-holomorphic sequence of sections of
 the vector bundles $\C^N \ox E \otimes L^{\otimes k}$ and $\a>0$. Then 
 there exists another sequence $\s_k$ verifying:
\begin{enumerate}
  \item $|s_{k}-\sigma_{k}|_{C^1,g_k}< \alpha$.
  \item $\sigma_{k}$ is $\gamma$-grassmannizable and $\gamma$-generic 
    of order $n$ for some $\gamma>0$.
\end{enumerate}
 Moreover, the result holds for one-parameter families of sections
 where the sections and the compatible almost complex structures
 depend continuously on $t\in [0,1]$.
\end{proposition}

\begin{lemma} \label{lem:emb->sec2}
  Let $\p_k:M\to \Gr(r,N)$ be a sequence of asymptotically 
  holomorphic embeddings with $\p_k^*\SU=E\ox L^{\ox k}$. 
  Then there exists a sequence of asymptotically holomorphic 
  sections $s_k$ of $\C^N\ox E \ox L^{\ox k}$, for $k$ large enough, 
  which is $\g$-grassmannizable and $\g$-generic of order $n$, for some
  $\g>0$, such that $\p_k=\Gr(s_k)$. The same holds for continuous
  one-parameter families of embeddings and compatible almost complex
  structures. 
\end{lemma}

{\bf Proof of Theorem \ref{existence21}.}
Note that the last property is obvious by the construction.
Let us begin with an asymptotically $J$-holomorphic sequence $\s_k$ of 
sections of the bundles $\C^N\otimes E\otimes L^{\otimes k}$ and perturb 
it using Proposition \ref{key2} to obtain an asymptotically holomorphic 
$\gamma$-grassmannizable and $\gamma$-generic of order $n$ sequence 
of sections $s_k$. The first property implies that $\phi_k=\Gr (s_k)$
is well-defined, the second that it is an immersion. To get an embedding
we use that $2\dim M< \dim \Gr(r,N)=2r(N-r)$  to find a generic 
$C^p$-perturbation of norm less than $O(k^{-1/2})$
to get rid of the self-intersections and keeping the asymptotic 
holomorphicity, the grassmannizability and the genericity of order $n$.
Now we only have to check that
the sequence $\phi_k=\Gr (s_k)$ verifies the required conditions
in Definition \ref{def:ah}. 

Choose a point $x\in M$ and trivialize $E$ in a neighborhood of $x$ by
fixing an orthonormal basis $e_1,\ldots, e_r$.
Now by a rotation with an element of 
$U(N)$ acting on $\C^N$ and an element of $U(r)$ acting
on $E$, we can assure that 
\begin{equation}
 s_k(x)= \left( \begin{array}{ccccccc}
 s_k^{11}(x) & 0 & \ldots & & & \ldots & 0 \\
 0 & s_k^{22}(x) & 0 & \ldots & & & 0 \\
 0 & \ldots & \ddots & 0 & \ldots & & 0 \\
 0 & \ldots & & s_k^{rr}(x) & 0 & \ldots & 0 \end{array} \right) 
\label{matriz2}
\end{equation}
where $s_k^{ij}$ are sections of $L^{\ox k}$.
This corresponds to an isometric transformation
of $\Gr(r,N)$. The $\gamma$-grassmannizable property implies that 
$|s_k^{11}\cdots s_k^{rr}|\geq \gamma$. 
By the asymptotic holomorphicity bounds it is $|s_k|=O(1)$, so 
that $|s_k^{ii}| \geq \gamma/C$, for some universal constant $C$.
Therefore on a ball $B_{g_k}(x,c)$ of fixed universal radius $c$,
the first $r\x r$ minor of $s_k(y)$ has an inverse of norm bounded 
by $C' \gamma^{-1}$, for some universal constant $C'$.

Let $v_1,\ldots, v_N$ be the canonical basis of $\C^N$.
As $\p_k(x)=\Pi_0=\left[ \begin{array}{c} v_1 \\ \vdots \\ 
v_r\end{array}\right]$,
we consider the standard local chart 
for $\Gr(r,N)$ around $\Pi_0$ for the open set
$U_0=\{\Pi\big| \Pi \cap [v_{r+1},\ldots,v_N]=\{0\}\}$, given by
\begin{eqnarray*}
 \P_0 : U_0  &\to & \C^{r\x (N-r)} \\
   \left[ \begin{array}{ccc} s_{11} &\cdots &s_{1N} \\
  \vdots &\ddots &\vdots \\ s_{r1} &\cdots & s_{rN} \end{array} \right] 
   &\mapsto & \left( \begin{array}{ccc}
   s_{11} & \ldots & s_{1r} \\ \vdots   &  \ddots & \vdots  \\
   s_{r1} & \ldots & s_{rr}\end{array} \right)^{-1} \left( \begin{array}{cccc}
   s_{1,r+1} &s_{1,r+2} & \ldots & s_{1N} \\ \vdots  &  & \ddots & \vdots  \\
   s_{r,r+1} &s_{r,r+2}  & \ldots & s_{rN} \end{array} \right)
\end{eqnarray*}
It is easy to check that $\P_0$ is an isometry at the point $\Pi_0$.

The application $f_k= \P_0 \circ \p_k$ is given by
\begin{eqnarray*}
  f_k:  B_{g_k}(x,c) & \to & \C^{r\times (N-r)} \\
  y & \mapsto & \left( \begin{array}{ccc}
  s_k^{11}(y) & \ldots & s_k^{1r}(y) \\
  \vdots   & \ddots & \vdots  \\
  s_k^{r1}(y) & \ldots & s_k^{rr}(y) \end{array} \right)^{-1} 
  \left( \begin{array}{ccc}
  s_k^{1,r+1}(y)  & \ldots & s_k^{1N}(y) \\
  \vdots  &  \ddots & \vdots  \\
  s_k^{r,r+1}(y) & \ldots & s_k^{rN}(y) \end{array} \right)
\end{eqnarray*}

We can compute the bounds required in Definition \ref{def:ah} using $f_k$ 
instead of $\p_k$. Now the arguments in the proof of Theorem \ref{existence1}
carry over verbatim. For the isotopy result we use Lemma \ref{lem:emb->sec2}.
\hfill $\Box$

\subsection{Construction of $\gamma$-grassmannizable sections.}
 Our objective is to prove the following perturbation result:
\begin{proposition} \label{grassmannizable}
 Suppose $N>n+r-1$. 
 Let $s_k$ be an asymptotically $J$-holomorphic sequence of sections of the
 vector bundles $\C^N\ox E\otimes L^{\otimes k}$ which is
 $\g$-grassmannizable of order $l$, for some $\g>0$. Then given $\alpha>0$, 
 there exists an asymptotically 
 $J$-holomorphic sequence of sections $\sigma_k$ verifying:
\begin{enumerate}
  \item $|s_k-\sigma_k|_{C^1,g_k}< \alpha$.
  \item $\sigma_k$ is $\eta$-grassmannizable of order $l+1$ for some $\eta>0$.
\end{enumerate}
 Moreover, the result can be extended to continuous one-parameter families 
 depending continuously of $t\in [0,1]$.
\end{proposition}

{\bf Proof.}
 Again we use the globalization argument described in 
 Proposition \ref{globalizate}. Let us do the non-parametric case, the other
 one being a trivial extension by now.
 Define the local and $C^0$-open property $\CrP (\epsilon,x)$ as 
 $|\bigwedge^{l+1} s_k(x)|> \epsilon$. We only need to find for a point 
 $x\in M$ a section $\tau_{k,x}$ with Gaussian decay away from $x$, 
 assuring that $s_k+\tau_{k,x}$ verifies $\CrP(\eta,
 y)$ in a ball of universal $g_k$-radius $c$. 

 Choose a point $x\in M$. Fix an orthonormal basis $e_1,\ldots, e_r$
 trivializing $E$ in a neighbourhood of $x$, so $s_k$ may be interpreted
 as a morphism $\C^N \to \C^r \ox L^{\ox k}$.
 By a rotation with an element 
 of $U(N)$ on $\C^N$ and an element of $U(r)$ on $E$, we can assure that 
$$
 s_k(x)= \left( \begin{array}{ccccccc}
 s_k^{11}(x) & 0 & \ldots & & & \ldots & 0 \\
 0 & s_k^{22}(x) & 0 & \ldots & & & 0 \\
 0 & \ldots & \ddots & 0 & \ldots & & 0 \\
 0 & \ldots & & s_k^{rr}(x) & 0 & \ldots & 0 \end{array} \right)
$$
 with $|s_k^{11}(x)\cdots s_k^{ll}(x)|\geq \gamma$. So
 $|s_k^{11}\cdots s_k^{ll}|> \gamma/2$ on a ball $B_{g_k}(x,c)$ of fixed 
 radius $c$.
 Let $\sref$ be the sections given by Lemma \ref{localized} and
 define $\theta_k=s_k^{11}\cdots s_k^{ll} \sref$. Clearly
 $|\theta_k| > c_s\gamma/2$ on $B_{g_k}(x,c)$. Consider the family
 of functions
$$
  M_k^p= s_k^{11}\cdots s_k^{ll} s_k^{l+1,p}, \quad l+1 \leq p\leq N.
$$
 These are components of $\bigwedge^{l+1} s_k$. If we perturb $s_k$
 so that the norm of $M_k=(M_k^{l+1},\ldots, M_k^N)$ is bigger than
 $\eta=c' \d (\log (\d^{-1}))^{-p}$ then we have finished. For this
 we define $g_k=(g_k^{l+1},\ldots, g_k^N)=\left( 
 \frac{M_k^{l+1}}{\theta_k},\ldots, \frac{M_k^{N}}{\theta_k} \right)
 =\left( \frac{s_k^{l+1,l+1}}{\sref},\ldots,\frac{s_k^{l+1,N}}{\sref}\right)$.
 We obtain, scaling the coordinates by universal constants if necessary, 
 $g_k: B^+ \to \C^{N-l}$ which is asymptotically holomorphic. As 
 $n<N-l$, we can find $|w_k|<\d$ such that $|g_k-w_k|>\d 
 (\log(\d^{-1}))^{-p}$. Then we obtain that
 $|(M_k^{l+1}-w_k^{l+1}\theta_k,\ldots, M_k^N-w_k^N\theta_k)| >\eta=
 c'\d(\log (\d^{-1}))^{-p}$, for some universal $c'$. This
 perturbation term is achieved by adding the section
 $\tau_{k,x}= -(0,\stackrel{(l)}{\cdots}, 0, w_k^{l+1}e_{l+1} 
 \sref, \ldots, w_k^N e_{l+1} \sref)$ of the bundles 
 $\C^N\ox E\ox L^{\ox k}$. This finishes the proof. \hfill $\Box$

\begin{remark}
 We cannot improve the condition $N>n+r-1$ in 
 Proposition \ref{grassmannizable}. As we shall see in 
 Section \ref{determinantal}, we expect the locus of points of $M$ where the 
 rank of $s_k:\underline{\C}^N \to E\ox L^{\ox k}$ is not maximum to have 
 codimension $N-r+1$.
\end{remark}

\subsection{Inductive construction of sections $\gamma$-generic of order $l$}
Now we study the problem of perturbing the sequence $s_k$ to achieve 
genericity of order $n$. The result to be proved is the following.

\begin{proposition}
 Suppose $r(N-r)>2n$. 
 Let $s_k$ be an asymptotically $J$-holomorphic sequence of sections of 
 the vector bundles $\C^{N}\otimes E\ox L^{\otimes k}$, which is 
 $\g$-grassmannizable and $\g$-generic of order $l$. Then given $\alpha>0$, 
 there exists an asymptotically
 $J$-holomorphic sequence of sections $\sigma_k$ verifying:
\begin{enumerate}
  \item $|s_k-\sigma_k|_{C^1,g_k}< \alpha$.
  \item $\sigma_k$ is $\eta$-generic of order $l+1$ for some $\eta>0$.
\end{enumerate}
 Moreover, this result can be extended to continuous one-parameter families of
 sections and almost complex structures.
\end{proposition}

{\bf Proof.}
Define the property $\CrP (\epsilon, x)$ for a section $s_k$
which is $\g/2$-grassmannizable and $\g/2$-generic of order $l$ as
$|\bigwedge^{l+1} \partial \Gr(s_k)(x)|>\epsilon$. A perturbation of 
our initial section verifies the hypothesis if we perturb by adding 
sections of $C^1$ norm smaller than $\g/2C$, $C$ some universal constant.
For applying Proposition \ref{globalizate} we need to build, for
$0<\d< \g/2Cc''$, a local perturbation $\tau_{k,x}$ with $|\tau_{k,x}|<c''\d$
and Gaussian decay with the property $\CrP(\eta, y)$ on $B_{g_k}(x,c)$
with $\eta=c'\d(\log(\d^{-1}))^{-p}$. 

Choose a point $x\in M$. By a rotation with an element of 
$U(N)$ acting on $\C^N$ and an element of $U(r)$ acting
on $E$, we can assure that $s_k(x)$ is as in \eqref{matriz2}.
By the $\g$-grassmannizability, $|s_k^{11}(x)\cdots s_k^{rr}(x)| \geq \g$.
The asymptotically holomorphic bounds imply that $|s_k|=O(1)$, so 
that $|s_k^{ii}(x)| \geq \gamma/C$ for some universal constant $C$.
There is a fixed universal radius $c$ such that
the first $r\x r$ minor of $s_k(y)$ has an inverse of norm bounded 
by $C' \gamma^{-1}$, for some universal constant $C'$, on $B_{g_k}(x,c)$.
Then we can use the trivialization $\P_0$ to define the applications
\begin{eqnarray*}
  f_k:  B_{g_k}(x,c) & \to & \C^{r\times (N-r)} \\
  y & \mapsto & \left( \begin{array}{ccc}
  s_k^{11}(y) & \ldots & s_k^{1r}(y) \\
  \vdots   &  & \vdots  \\
  s_k^{r1}(y) & \ldots & s_k^{rr}(y) \end{array} \right)^{-1} 
  \left( \begin{array}{cccc}
  s_k^{1,r+1}(y) & \ldots & s_k^{1N}(y) \\
  \vdots  &  &  \vdots  \\
  s_k^{r,r+1}(y) & \ldots & s_k^{rN}(y) \end{array} \right)
\end{eqnarray*}
Now consider the sections $\sref$ of Lemma \ref{localized}.
We define the applications
\begin{eqnarray*}
  \tilde{f}_k:  B_{g_k}(x,c) & \to & \C^{r\times (N-r)} \\
  y & \mapsto &   \frac{1}{\sref(y)} \left( \begin{array}{cccc}
  s_k^{1,r+1}(y) &s_k^{1,r+2}(y) &  \ldots & s_k^{1N}(y) \\
  \vdots  &  &  \vdots  \\
  s_k^{r,r+1}(y) & s_k^{r,r+2}(y) & \ldots & s_k^{rN}(y) \end{array} \right)
\end{eqnarray*}
Clearly $f_k=\Psi \circ \tilde{f}_k$ where $\Psi:B_{g_k}(x,c) \to \GL(r,\C)$
satisfies
$|\Psi|=O(1)$, $|\Psi^{-1}|=O(1)$,
$|\nabla \Psi|=O(1)$ and $|\nabla \Psi^{-1}|=O(1)$. Therefore it is 
enough to get a perturbation which has
$|\bigwedge^{l+1} \bd \tilde{f}_k|>\eta$ on $B_{g_k}(x,c)$.

Spreading out the entries of the matrix $\tilde{f}_k$ in one row we can write
$\tilde{f}_k(y)=(\tilde{f}_k^{11}(y), \ldots, \tilde{f}_k^{r,N-r} (y))$. 
Using the local forms $dz^1_k,\ldots, dz^n_k$, we may write
$$ 
 \bd \tilde{f}_k = (u_k^{111}dz_k^1+u_k^{112}dz_k^2+\cdots+u_k^{11n}dz_k^{n}, 
 \,\ldots\, , u_k^{r,N-r,1}dz_k^1+\cdots + u_k^{r,N-r,n}dz_k^n ), 
$$
for some $u_k^{ijl}$. Using a unitary transformation of $U(n)$ 
on the complex Darboux coordinate chart and relabeling horizontally the
coordinates, we can suppose that
\begin{equation} 
 \bd \tilde{f}_k(x)= \left( \begin{array}{ccccccc}
  u_k^{11}(x) & * & \ldots & & & \ldots & * \\
  0 & u_k^{22}(x) & * & \ldots & & & * \\
  0 & \ldots & \ddots & * & \ldots & & * \\
  0 & \ldots & 0 & u_k^{nn}(x) & * & \ldots & * \end{array} \right), 
\label{matriz3}
\end{equation}
 where $|u_k^{11}(x)\cdots u_k^{ll}(x)|>\g/C_0$, $C_0$ a universal constant.
 The relabeling is given by a function $\a\in\{1,\ldots, r(N-r)\} \mapsto
 (i(\a),j(\a)) \in \{1,\ldots,r\}\x\{1,\ldots, N-r\}$.
 Shrinking $c$ if necessary we can assure that $|(\partial \tilde{f}_k^1\wedge
 \cdots \wedge \partial \tilde{f}_k^l)_{dz_k^1\wedge \cdots 
 \wedge dz_k^l}|> \g/2C_0$ for all the points of the 
 ball $B_{g_k}(x,c)$. Now we construct the $(l+1)$-form
$$ 
  \theta_k(y)= (\partial \tilde{f}_k^1\wedge \cdots \wedge \partial 
  \tilde{f}_k^l)_{dz_k^1\wedge \cdots \wedge dz_k^l} \wedge dz_k^{l+1}. 
$$
and the family of $(l+1)$-forms
$$
  M_k^p= (\partial \tilde{f}_k^1\wedge \cdots \wedge \partial
  \tilde{f}_k^l \wedge  \partial \tilde{f}_k^p)_{dz^1_k \wedge 
  \cdots \wedge dz_k^l \wedge dz_k^{l+1}}, 
  \quad l+1\leq p \leq r(N-r),
$$
 which are components of $\bigwedge^{l+1} \bd \tilde{f}_k$. If we 
 perturb so that the norm of $M_k=(M_k^{l+1}, \ldots , M_k^{r(N-r)})$ gets
 bigger than $\eta=c'\d (\log (\d^{-1}))^{-p}$ then we are done.
 We define the following sequence of asymptotically holomorphic
 applications: $h_k=\left(
 \frac{M_k^{l+1}}{\theta_k}, \ldots, \frac{M_k^{r(N-r)}}{\theta_k}\right)$.
 So we obtain, scaling the coordinates by universal constants if necessary,
 $h_k: B^+ \to \C^{r(N-r)-l}$ which is asymptotically holomorphic.
 As $n<r(N-r)-l$ we can find $|w_k|<\delta$ such that 
 $|h_k-w_k|>\delta (\log (\delta ^{-1}))^{-p}$.
 Thus $|(M_k^{l+1}-w_k^{l+1}\theta_k,
 \ldots, M_k^{r(N-r)}-w_k^{r(N-r)}\theta_k)|> \eta=c'\d(\log(\d^{-1}))^{-p}$. 
 The perturbation term $-(w_k^{l+1} \theta_k, \ldots, w_k^{r(N-r)} \theta_k)$
 is achieved by adding the section
 $$
  \tau_{k,x}= -(0,\stackrel{(r)}{\ldots}, 0, 
  \hspace{-3mm} \sum_{j(\a)=r+1,\a>l} \hspace{-3mm}w_k^{\a}z_{l+1}e_{i(\a)} 
  \sref, \ldots,\hspace{-3mm}\sum_{j(\a)=N,\a>l}\hspace{-3mm}
  w_k^{\a} z_{l+1} e_{i(\a)} \sref).
 $$ 
 This finishes the proof in the non-parametric case. The other case is left
 to the reader.
\hfill $\Box$

\subsection{Lifting asymptotically holomorphic embeddings in grassmannians}
\label{lifting2}
This Subsection is devoted to a proof of Lemma \ref{lem:emb->sec2}, 
which states that any asymptotically holomorphic embedding into a
grassmannian is of the form provided by Theorem \ref{existence21}.

{\bf Proof of Lemma \ref{lem:emb->sec2}.}
Suppose that we have a sequence of $\g$-asymptotically holomorphic
embeddings $\p_k:M\to \Gr(r,N)$, for some $\g>0$, with
$\p_k^*\SU=E\ox L^{\ox k}$, where $\SU$ is the universal rank $r$ 
bundle over the grassmannian. The dual of $\SU$ is given by
$$
  \SU^*=\{(\Pi,s) \big| s \in \Pi \} \subset \Gr(r,N) \x \C^N=
  \underline{\C}^N,
$$
interpreted as a sub-bundle of the trivial bundle $\underline{\C}^N$.
We consider the sequence of hermitian bundles, $E_k=\p_k^* \SU^* \ox E \ox
L^{\ox k}= \End E \subset \C^N\ox E \ox L^{\ox k}$.
We look for sequences of sections $s_k$ of $E_k$ which are 
$\s$-grassmannizable of order $n$ such that they are asymptotically 
holomorphic when considered as sections of $\C^N\ox E \ox L^{\ox k}$.
Let $S^l_k=\Tr (\bigwedge^l s_k)$, which is an asymptotically 
holomorphic sequence of sections of the trivial
vector bundle $\underline{\C}$. We want to prove that $|S^r_k|\geq \s$
for $k$ large. We shall prove that we can find sequences $s_k$ with
$|S^l_k|\geq \eta_l$, for some $\eta_l>0$, by induction on $l$.

Suppose that $s_k$ is an asymptotically holomorphic sequence
of sections of $E_k$ such that $|S^l_k|\geq \g$. 
Let $\CrP(\f,x)$ be the $C^1$-open property for sequences of 
sections $s_k$ of $E_k$ given as $S^{l+1}_k=\Tr 
(\bigwedge^{l+1} s_k)$ is $\f$-transverse to $0$ at $x$.

Let $x\in M$. We want to find a local perturbation with Gaussian
decay obtaining the property $\CrP(\eta,y)$ in a ball of
universal $g_k$-radius $c$ around $x$. For this, 
define the local sections $\s_k$ of 
$\p_k^*\SU^* \ox E \subset {\C}^N \ox E$ as follows. Locally,
$\s_k$ is a map
$$
  \s_k: B_{g_k}(x,c) \to \C^N \ox \C^r,
$$
such that for $y\in B_{g_k}(x,c)$, $\s_k(y)$ is a $N\x r$ matrix, i.e.\
a linear map $\s_k(y):\C^N \to \C^r$. The point $\p_k(y)\in \Gr(r,N)$ 
corresponds to the image of the embedding $\s_k(y)^T:\C^r \to \C^N$. 
Note that one may identify the tangent space $T_{\p_k(y)}\Gr(r,N)$ 
to the set of linear maps $\C^N \to \C^r$ which are zero on 
$\im (\s_k(y)^T)$, i.e.\ maps $\varphi$ such that $\varphi\, \s_k(y)^*=0$. 
With this, $\nabla \s_k=\nabla \p_k+ (\nabla \s_k \s_k^*) \s_k$.
So it is natural to require $(\nabla_r \s_k(y))\s_k(y)^*=0$, for any
$y\in B_{g_k}(x,c)$, where $r$ is the radial vector field from $x$.
We fix $\s_k(x)$ satisfying $\s_k(x)\s_k(x)^*=\id$.
This determines $\s_k$ uniquely.
The following bounds are proved as in Subsection \ref{lifting},
 \begin{eqnarray*}
  & & \s_k(y)\s_k(y)^*=\id, ~~~~ |\s_k(y)|=O(1), ~~~~ 
  |\nabla \s_k(y)|=O(1+d_k(x,y)), \\
  & & |\bbd \s_k(y)|=O(k^{-1/2}(1+d_k(x,y))), ~~~~ 
  |\nabla \bbd \s_k(y)|=O(k^{-1/2}(1+d_k(x,y))). 
 \end{eqnarray*}

Trivialize $E$ in a ball $B_{g_k}(x,c)$, so that $s_k/\sref$ can
be considered as an application $B_{g_k}(x,c) \to \C^{r\x N}$.
Define the application
$$
  f_k= \frac{s_k \s_k^*}{\sref}: B_{g_k}(x,c) \to \C^{r\x r},
$$
so that $f_k\,\s_k=s_k/\sref$. Then $f_k$
is asymptotically holomorphic and we may check
property $\CrP(\eta,y)$ for $f_k$ instead of $s_k$. 
Let $F_i=\Tr (\bigwedge^i f_k)$, so that $|F_l| \geq C\g$ for
some universal constant $C$. For any $w\in \C$ we have
$$
  \Tr (\bigwedge\nolimits^{l+1} (f_k-w\id))= F_{l+1}-w(n-l) F_l+ 
  w^2\binom{n-l+1}{2} F_{l-1} +\ldots +(-w)^{l+1} \binom{n}{l+1}F_0
$$
By the standard argument, we may obtain $|w|<\d$ such that 
$\frac{F_{l+1}}{F_l}-w$ is $\eta$-transverse to $0$, with 
$\eta=\d(\log (\d^{-1}))^{-p}$, in $B_{g_k}(x,c)$.
Then it is easy to see that $\Tr (\bigwedge^{l+1} (f_k-\frac{w}{n-l}
\id)$ is $c'\eta$-transverse to $0$, where $c'$ is a universal
constant. The perturbation
$$
  \tau_{k,x}=-\frac{w}{n-l}\s_k \sref
$$
is a sequence of sections of $E_k=\p_k^*\SU^* \ox E\ox L^{\ox k}$, 
with Gaussian decay such that $|\tau_{k,x}|<c''\d$ and
$s_k+\tau_{k,x}$ satisfies $\CrP(\eta,y)$ for $y\in B_{g_k}(x,c)$,
with $\eta=c'\d(\log (\d^{-1}))^{-p}$. By Proposition \ref{globalizate}, 
there exists an asymptotically holomorphic sequence of sections
of $E_k$, which we denote by $s_k$ again, such that 
$S^{l+1}_k=\Tr(\bigwedge^{l+1} s_k)$ is $\eta$-transverse to $0$, 
for some $\eta>0$. For $k$ large enough, the zeroes of $S_k^{l+1}$ is 
a symplectic submanifold representing the trivial homology class, 
hence the empty set. So $|S_k^{l+1}|\geq \eta$.
This completes the proof. The extension to the one-parameter case is trivial.
\hfill $\Box$

\subsection{Zero sets of vector bundles}
Following the ideas of Subsection \ref{esti_inter} and using 
Proposition \ref{angle} we can prove the following two results

\begin{theorem} \label{good_inter2}
 Given $\p_k= \Gr (s_k)$, where $s_k$ is a sequence of asymptotically 
 holomorphic sections of $\C^N \ox E\ox L^{\ox k}$, which are 
 $\g$-grassmannizable and $\g$-generic of order $n$, for some $\g>0$.
 Fix a holomorphic submanifold $V$ in $\Gr(r,N)$. Then for any $\a>0$ 
 there exists a sequence of asymptotically holomorphic sections 
 $\s_k$ of $\C^N \otimes E \otimes L^{\otimes k}$ such that
  \begin{enumerate}
    \item $|\s_k-s_k|_{g_k,C^1}<\a$.
    \item $\Gr (\s_k)$ is an $\eta$-asymptotically holomorphic embedding
    in $\Gr(r,N)$ which is $\epsilon$-transverse to $V$, with $\eta>0$ and
    $\epsilon>0$  independent of $k$. In the case $\dim M + \dim V < 2r(N-r)$ 
     we have that $d_{Gr}(\phi_k(M),V))>\epsilon$, for $k$ large enough. 
  \end{enumerate}
 Moreover the result can be extended to one-parameter continuous families
 of complex submanifolds $(V_t)_{t\in [0,1]}$, taking in this case as 
 starting point a continuous family $\phi_{t,k}=\Gr(s_{t,k})$, where
 $s_{t,k}$ is a continuous family of $\gamma$-grassmannizable and $\g$-generic
 of order $n$ sequences asymptotically $J_t$-holomorphic.
\end{theorem}

{\bf Proof.}
The proof is similar to that of Theorem \ref{good_inter1}. We just briefly
point out the differences. For simplicity we suppose that the codimension 
of $V$ is $1$. 

For $x\in M$, we may suppose that $s_k(x)$ is as in \eqref{matriz2}.
We use the chart $\P_0$ to get the local maps 
\begin{eqnarray*}
  f_k= \P_0\circ \p_k:  B_{g_k}(x,c) & \to & \C^{r\times (N-r)} \\
  y & \mapsto & \left( \begin{array}{ccc}
  s_k^{11}(y) & \ldots & s_k^{1r}(y) \\
  \vdots   &  & \vdots  \\
  s_k^{r1}(y) & \ldots & s_k^{rr}(y) \end{array} \right)^{-1} 
  \left( \begin{array}{cccc}
  s_k^{1,r+1}(y) & \ldots & s_k^{1N}(y) \\
  \vdots  &  &  \vdots  \\
  s_k^{r,r+1}(y) & \ldots & s_k^{rN}(y) \end{array} \right)
\end{eqnarray*}
This time we have a vector $v \in \C^{r\x (N-r)}$.
We define the functions $h_k:B_{g_k}(x,c) \to \C$ by the condition
$$
  f \left( f_k +r_k \sref \left( \begin{array}{ccc}
  s^{11}_k &\cdots &  s^{1r}_k  \\
  \vdots & \ddots & \vdots  \\
  s^{r1}_k & \cdots & s^{rr}_k \end{array} \right)^{-1} 
  \left( \begin{array}{ccc}
  s^{11}_k(x) &\cdots &  0  \\
  \vdots & \ddots & \vdots  \\
  0 & \cdots & s^{rr}_k(x) \end{array} \right) 
  \left( \begin{array}{ccc}
  v^{11}&\cdots &  v^{1,N-r}  \\
  \vdots &\ddots & \vdots  \\
  v^{r1} & \cdots &  v^{r,N-r} \end{array} \right)\right)=0,
$$
and prove that they are asymptotically holomorphic.
Then we find $|w_k|<\d$ such that $h_k-w_k$ is $\eta$-transverse
to $0$ with $\eta=c'\d(\log(\d^{-1}))^{-p}$. 
Finally the perturbation will be
$$
 \tau_{k,x} = - \left( \begin{array}{cccccc}
  0 & \cdots & 0 & w_k v^{11}\sref & \cdots & w_k v^{1,N-r}\sref  \\
  \vdots & \ddots & \vdots & \vdots & \ddots & \vdots \\
  0 & \cdots & 0 & w_k v^{r1}\sref & \cdots & w_k v^{r,N-r} \sref 
 \end{array} \right).
$$
The arguments run parallel to those in the
proof of Theorem \ref{good_inter1}, although the 
constants have to be arranged suitably, but we leave this task to
the careful reader. 
\hfill $\Box$

We call universal planes to the zero sets of algebraic sections 
transverse to zero of the universal bundle $\SU$
over the grassmannian $\Gr(r,N)$. Now we can deduce the main result of 
\cite{Au97}.

\begin{corollary}
 Let $(M, \omega)$ be a compact symplectic manifold of integer class.
 Let $E$ be a hermitian vector bundle over $M$. Then for $k$ large 
 enough there exist symplectic submanifolds obtained as zero sets 
 of the bundles $E\ox L^{\ox k}$. Moreover, perhaps by increasing $k$,
 we can assure that all the symplectic submanifolds constructed
 as transverse intersections of asymptotically holomorphic
 sequences with a fixed universal plane are isotopic.
 The isotopy can be made through symplectomorphisms.
\end{corollary}
 
The proof follows the steps of the proof of Corollary \ref{bertini}. 
Remark also that the result is a corollary of Theorem \ref{main_depend}
to be proved in Section \ref{determinantal}.

\section{Determinantal submanifolds of closed symplectic manifolds} 
\label{determinantal}

Let $(M,\o)$ be a symplectic $4$-manifold of integer class, endowed with
a compatible almost complex structure. Let $E$ and $F$ be two vector bundles of
ranks $r_e$ and $r_f$, respectively. Recall that for any morphism $\varphi:
E\to F$ we have defined in Definition \ref{defn:determ}
the $r$-determinantal set as 
$$
 \Sigma^r(\varphi)= \{ x \in M \big| \rank\, \varphi_x=r \}. 
$$
We want to prove Theorem \ref{main_deter},
which allows to construct $\Sigma^r(\varphi)$ as a symplectic submanifold, 
after twisting $E$ and $F$ with large powers of $L$.
The solution to this problem goes through embedding $M$ in a product 
of two grassmannians and cutting its image with suitable ``generalized
Schur cycles''. We shall do this in next Section.

\begin{remark}
 A direct approach to proving Theorem \ref{main_deter} 
 consists on reducing it to Auroux' case by taking 
 the $r$-fold wedge product of $\varphi_k$,
\begin{eqnarray*}
  \bigwedge\nolimits^r \varphi_k: \bigwedge\nolimits^r E\ox (L^*)^{\ox k}  & 
  \to & \bigwedge\nolimits^r F\ox L^{\ox k} \\
   s_1 \wedge \cdots \wedge s_r & \mapsto & 
   \varphi_k(s_1)\wedge \cdots \wedge \varphi_k(s_r).
\end{eqnarray*}
 So the zero set of $\bigwedge^r \varphi_k$ is generically a stratified 
 submanifold $\Sigma^0(\varphi_k) \bigcup \ldots \bigcup \Sigma^r(\varphi_k)$.
 If we suppose that $\varphi_k$ is an asymptotically $J$-holomorphic sequence
 of sections of the bundle $E^* \ox F \ox L^{\ox 2k}$, one could try to use
 Donaldson's techniques to obtain a new sequence of sections transverse in an
 adequate sense to assure the symplecticity. The following example shows the
 main obstacle to this approach. Take a symplectic $4$-manifold in the
 hypothesis of Theorem \ref{main_deter} with two hermitian vector 
 bundles $E$ and $F$ of rank $2$. Using Auroux' techniques we can assure that
 the zero sets of $\varphi_k$ are $\eta$-transverse to $0$, for some $\eta>0$.
 When we go to $\bigwedge^2 \varphi_k$, the condition to be satisfied is
$$ 
 |\bbd \bigwedge\nolimits^2 \varphi_k| < |\bd \bigwedge\nolimits^2 \varphi_k |.
$$
 However, at any $x \in Z(\varphi_k)$ we obtain $|\bbd \bigwedge^2 \varphi_k
 (x)| =|\bd \bigwedge^2 \varphi_k(x)|=0$, so we cannot impose a global
 transversality property for the section $\bigwedge^2 \varphi_k$. This case
 is very similar to that in \cite{Do99} and can be treated with an ``ad hoc''
 argument, but more general cases do not admit a treatment based on the use
 of {\it normal forms} of the singularities,
 because for higher dimensions the problem becomes intractable \cite{Ar82}.
\end{remark}

\subsection{Bigrassmannian embeddings} \label{bigr}
 The idea to prove Theorem \ref{main_deter} is based in the following 
 observations. Choose two sequences of sections $s_k^e$ and $s_k^f$ of
 the bundles $\C^N\otimes E^*\otimes L^{\otimes k}$ and 
 $\C^N\otimes F \otimes L^{\otimes k}$ respectively, which are 
 $\gamma$-grassmannizable and $\gamma$-generic of order $n$, for some $\g>0$,
 providing by Theorem \ref{existence21}, asymptotically holomorphic sequences of
 embeddings $\Gr(s_k^e)$ and $\Gr(s_k^f)$ of $M$ in $\Gr(r_e,N)$ and
 $\Gr(r_f,N)$, respectively, for $N$ a large integer number.

 Performing the cartesian product we obtain an asymptotically
 holomorphic sequence of embeddings of $M$ into the bigrassmannian
 $\Bi(r_e,r_f,N)=\Gr(r_e,N)\x\Gr(r_f,N)$,
 $$
  \p_k=\Gr(s_k^e)\x\Gr(s_k^f):M\to\Gr(r_e,N)\x\Gr(r_f,N)=\Bi(r_e,r_f,N).
 $$
 Let $\SU_e$ and $\SU_f$ be the universal bundles over $\Gr(r_e, N)$ and
 $\Gr(r_f,N)$ respectively, which are very ample. Define 
 $\pi_e: \Bi(r_e,r_f,N) \to \Gr(r_e,N)$ as the projection onto the 
 first factor (and analogously $\pi_f$). Therefore $\SU_{ef}=\pi_e^*(\SU_e) 
 \otimes \pi_f^*(\SU_f)=\SU_{ef}$ is very ample on $\Bi(r_e,r_f,N)$.
 Recall that $\Gr(s_k^e)^* (\SU_e)= E^* \otimes L^{\otimes k}$ and 
 $\Gr(s_k^f)^* (\SU_f) = F \otimes L^{\otimes k}$. Then $\phi_k^* 
 \SU_{ef}= E^* \otimes F \ox L^{\ox 2k}$.
 $\SU_{ef}$ has a holomorphic section $s$ verifying that: 
\begin{enumerate}
  \item $D_r=\Sigma^r(s)$ is an open complex submanifold in $\Bi(r_e,r_f,N)$.
  \item $\codimC D_r= (r_e-r)(r_f-r)$. 
\end{enumerate}

 If we assure that, for each $r$, $\phi_k$ is transverse to
 $D_r$ with an angle $\epsilon>0$ independent of $k$, we have 
 finished the proof of Theorem \ref{main_deter} by 
 Proposition \ref{angle}. This is carried out as follows.

\begin{lemma} \label{strata_bor}
 Let $\phi_k: M \to \Bi(r_e, r_f, N)$ be a $\gamma$-asymptotically 
 holomorphic sequence of embeddings. Suppose that $\phi_k$ is 
 $\s$-transverse to $D_{r}$. Then there exists $\epsilon>0$, 
 depending only on $\gamma$, $\s$ and the universal bounds of the 
 derivatives of the sequence, such that $\phi_k$ is $\s/2$-transverse to
 $D_{r'}$, $r'>r$, when we restrict to an $\epsilon$-neighborhood of $D_{r}$.
\end{lemma}

In other words we do not have to care about the behaviour of the angle 
near the border of the strata.

{\bf Proof.}
Choose a point $x\in D_{r}\bigcap \phi_k(M)$. Recall that by 
$\s$-transversality, the minimum angle between $T_xD_r$ and $T_x\phi_k(M)$ 
is greater than $\s$.
We trivialize $\Bi(r_e, r_f, N)$ by a chart $\Phi_0$ defined as the 
cartesian product of two standard charts in the grassmannians, which 
is an isometry at the origin and verifies 
that $\Phi_0(x)=0$, namely,
$$ 
 \Phi_0: \Bi(r_e, r_f, N) \to \C^{r_e(N-r_e)}\times \C^{r_f(N-r_f)}. 
$$
Since $D_r$ is contained in the closure of $D_{r'}$, we have
\begin{equation}
  |y|< \delta \Rightarrow \angle_M(T_0\Phi_0(D_r), T_y\Phi_0(D_{r'}))< 
  c_D\delta,~~~~ \forall y \in B(0, c_u)\bigcap \Phi_0(D_{r'}). \label{prox1}
\end{equation}
The angles are measured with respect to the standard euclidean metric which is
close to that induced by the bigrassmannian if we choose $c_u$ small 
enough. Here $c_D$ is universal. Also by the asymptotic 
holomorphicity bounds of $\p_k$ we know that
\begin{eqnarray}
 & & |y|<\delta \Rightarrow \angle_M(T_0\Phi_0(\phi_k(M)), 
   T_y\Phi_0(\phi_k(M)))< c_{\phi}\delta, \nonumber \\ 
 & & \forall y \in B(0, c_u)\bigcap \Phi_0(\p_k(M)), \label{prox2}
\end{eqnarray}
where $c_{\phi}$ is universal. 
Now Proposition \ref{vari_min} says that
$$
  \angle_m(T_0\P_0(D_r),T_0\P_0(\p_k(M))) \leq  
   \angle_M(T_0\P_0(D_r),T_y\P_0(D_{r'}))+ \qquad 
$$
$$ 
  \qquad + \angle_m(T_y\P_0(D_{r'}),T_y\P_0(\p_k(M)))+ 
   \angle_M(T_y\P_0(\p_k(M)),T_0\P_0(\p_k(M))). 
$$
Using inequalities \eqref{prox1} and \eqref{prox2} and
remembering that all the angles have to be 
measured with respect to the bigrassmannian metric (which
is related to the standard metric in the ball $B(0,c_u)$
by non zero universal constants), we get the required result.
\hfill $\Box$

With Lemma \ref{strata_bor} the proof of Theorem \ref{main_deter} reduces 
to the following result, whose proof is similar to that of
Theorem \ref{good_inter2}.

\begin{proposition} \label{strata_int}
 Let $s_k^e$ and $s_k^f$ be two asymptotically holomorphic sequences of the 
 vector bundles $E^*\ox L^{\ox k}$ and $F\ox L^{\ox k}$ which are 
 $\gamma$-grassmannizable and $\gamma$-generic of order $n$, defining
 so an asymptotically holomorphic embedding in $\Bi(r_e, r_f,N)$. Fix an 
 algebraic open 
 submanifold $V$ in $\Bi(r_e, r_f,N)$ with compactification $\bar{V}=V 
 \bigcup W$. Then for any $\f, \a>0$, there exist $\eta>0$ and 
 two asymptotically holomorphic sequences $\sigma_k^e$ and $\sigma_k^f$ 
 of sections of the vector bundles $E^*\ox L^{\ox k}$ and $F\ox L^{\ox k}$ 
 respectively, verifying:
\begin{enumerate}
 \item $|\sigma_k^e-s_k^e|_{g_k,C^1}<\a$ and 
  $|\sigma_k^f-s_k^f|_{g_k,C^1}<\a$.
 \item $\phi_k=\Gr(\sigma_k^e)\times\Gr(\sigma_k^f)$ is a sequence of 
  $\eta$-asymptotically holomorphic embeddings in $\Bi(r_e, r_f, N)$.
 \item Denoting by $V_{\epsilon^-}$ the compact submanifold of $V$ 
  obtained by removing an $\epsilon$-neighborhood of $W$, we obtain 
  that $\phi_k$ is $\eta$-transverse to $V_{\epsilon^-}$.
\end{enumerate}
 Moreover the result can be extended to continuous one-parameter families 
 of sections $(s_{k,t}^e)_{t\in [0,1]}$ and
 $(s_{k,t}^f)_{t\in [0,1]}$ providing embeddings to the bigrassmannian and 
 to continuous one-parameter families of open submanifolds $V_t$. Thus we 
 obtain continuous families of sequences $\sigma_{k,t}^e$ and 
 $\sigma_{k,t}^f$ verifying the required conditions.
\hfill $\Box$
\end{proposition}

\subsection{Dependence loci of sections of a vector bundle}
 Suppose that $E$ is an hermitian vector bundle of rank $n$ and consider
 $s_1,\ldots, s_m$ sections of $E$. Then we can interpret $s=(s_1,\ldots,s_m)$
 as a morphism of bundles $s: \underline{\C}^m \to E$. The $r$-determinantal
 set of $s$ is
$$
 \Sigma^r(s)= \{ x \in M \big| \dim [s_1(x),\ldots, s_m(x)]=r \}, 
$$
 and it is called the $r$-dependence locus of the sections
 $s_1,\ldots, s_m$.

\begin{theorem} \label{main_depend}
 Let $(M, \omega)$ be a closed symplectic manifold of integer class and
 let $E$ be a rank $n$ hermitian vector bundle. 
 Then, for $k$ large enough, there exist $s_k=(s_k^1,\ldots, s_k^m)$
 sections of $\C^m \ox E$ such that
\begin{enumerate}
 \item $\Sigma^r(s_k)$ is an open symplectic submanifold of $M$.
 \item $\codim \Sigma^r(s_k)= 2(m-r)(n-r)$. The set of manifolds 
  $\{\Sigma^r(s_k)\}_r$ constitutes a stratified submanifold.
\end{enumerate}
 Moreover, any two stratified submanifolds constructed 
 by the process in the proof below are isotopic.
\end{theorem}

{\bf Proof.}
 The proof is similar to the arguments developed in Subsection \ref{bigr}.
 Let $\SU$ be the universal bundle over $\Gr(n, N)$ and
 consider $m$ holomorphic sections $s_1,\ldots, s_m$ verifying that: 
\begin{enumerate}
  \item $D_r=\Sigma^r(s)$ is an open complex submanifold in $\Gr(n,N)$.
  \item $\codimC D_r= (m-r)(n-r)$. 
\end{enumerate}
 Now we choose a sequence of asymptotically holomorphic embeddings 
 $\p_k:M \to \Gr(n,N)$ such that $\p_k^*\SU=E\ox L^{\ox k}$.
 If we assure that, for each $r$, $\phi_k$ is transverse to
 $D_r$ with an angle $\epsilon>0$ independent of $k$, we have 
 finished the proof because of Proposition \ref{angle}. But we 
 may perturb $\p_k$ by using analogues of Lemma \ref{strata_bor} 
 and Proposition \ref{strata_int} for the case of just one grassmannian.
\hfill $\Box$

\subsection{Homology and homotopy groups of the determinantal submanifolds}
In this Subsection we prove a result concerning the topology of smooth determinantal 
submanifolds analogous to Proposition 39 in \cite{Do96} (symplectic Lefschetz hyperplane
theorem) and Proposition 2 in \cite{Au97}. The main result is

\begin{proposition} \label{lefs}
Let $E,F$ be vector bundles of ranks $r_e$, $r_f$, respectively, over a
closed symplectic manifold $(M, \o)$ of integer class and
let $D_r^k$ be a sequence of determinantal submanifolds constructed, by using the
vector bundles $E\otimes (L^*)^{\otimes k}$ and $F\otimes L^{\otimes k}$,
as a transverse intersection of an asymptotically holomorphic sequence of embeddings
in $\Bi(r_e, r_f,N)$ with the determinantal varieties of a fixed generic section $s$
of the universal bundle $\SU_{ef}$ over $\Bi(r_e,r_f,N)$. Assume
that the stratified determinantal submanifold has only one stratum $D_r^k$.
Then the inclusion $i: D_r^k\to M$ induces, for $k$ large enough, an isomorphism 
on homotopy groups $\pi_p$ for 
$p< \frac12 \dim D_r^k$ and a surjection on $\pi_p$ for $p=\frac12 \dim D_r^k$. 
The same property also holds for homology groups.
\end{proposition}

Remark that the asumption of only one stratum implies that $r=\min \{ r_e, r_f \}-1$ and
$2(r_e-r+1)(r_f-r+1)=4(|r_e-r_f|+2)>\dim M$. Along the proofs we will suppose that 
$r_e\geq r_f$, leaving the details of the other case to the reader. We proceed in several 
steps.

\subsubsection{Determinant vector spaces}
Let $V,W$ be vector spaces of dimensions $m$ and $n$ ($m\geq n$) respectively.
We need some results about the behaviour of the determinant vector space 
$\bigwedge^r (V^*)\otimes \bigwedge^r W$ associated to the vector space of 
linear morphisms $V^*\otimes W$. We define the $r$-fold wedge product 
$\bigwedge^r$ of a linear application $\varphi\in \Hom(V,W)$
as the linear application
\begin{eqnarray*}
\bigwedge\nolimits^r \varphi: \bigwedge\nolimits^r V & \to & \bigwedge\nolimits^r W \\
v_1 \wedge \cdots \wedge v_r & \to & \varphi(v_1)\wedge \cdots \wedge \varphi(v_r).
\end{eqnarray*}
Thus we obtain a non-linear map $\bigwedge^r: \Hom(V,W) \to 
\Hom(\bigwedge^r V, \bigwedge^r W)$. The previous definition extends in an obvious 
way to any pair of vector bundles $E$ and $F$ 
providing a non-linear map of vector bundles 
$\bigwedge^r: \Hom(E,F) \to \Hom(\bigwedge^r E, \bigwedge^r F)$.
With this notation a rank $r-1$ determinantal submanifold $D_{r-1}$ 
associated to a morphism $\varphi$ between vector bundles $E$ and $F$ is the set 
$$ 
  D_{r-1}= \{ x\in M: \bigwedge\nolimits^r \varphi(x)= 0 \}. 
$$

\begin{lemma}
Let $V,W$ be vector spaces of dimensions $m$ and $n$ ($m\geq n$) respectively, then the set 
$R(V,W)=\bigwedge^n(\Hom(V,W))-\{ 0 \}$ is a smooth open complex submanifold of 
$\Hom(\bigwedge^n V, \bigwedge^n W)$ of dimension $m-n+1$. Moreover 
$R(V,W)$ is invariant under multiplication by non-zero complex scalars, and so
given any point $d\in R(V,W)$ then, using the standard identification between a
vector space and its tangent space at a point, $d\in T_d R(V,W)$.
\end{lemma}

{\bf Proof.}
The last statement is obvious. For the first one, fix basis $(e_1, \ldots, e_m)$ in
$V$ and $(f_1, \ldots, f_n)$ in $W$. First notice that
$R(V,W)$ is invariant under the actions of the groups $\GL(V)$ and $\GL(W)$. Thus
for computing $T_{\wedge^n (\varphi)} R(V,W)$ we can restrict our atention to the point
\begin{equation}
\varphi= \sum_{i=1}^n f_i\otimes e_i^*. \label{chula}
\end{equation}
Notice that this is possible since the condition $\bigwedge^n \varphi\neq
0$ implies that the linear map $\varphi$ has rank $n$ and therefore suitable changes of basis 
provide the expression \eqref{chula}. Now, we only have to compute the images of the tangent
basis $\varphi_{ij}=\frac{d}{dt}|_{t=0}(\varphi+tb_{ij})$, 
where $b_{ij}=f_j\otimes e_i^*$. First 
assume that $i\leq n$, then we obtain
$$ 
 (\bigwedge\nolimits^n)_* \varphi_{ij} = \left\{ \begin{array}{ll} \varphi, & i=j, \\
   0, & i\neq j \end{array} \right. 
$$
However for the cases $i>n$ we obtain
$$ 
  (\bigwedge\nolimits^n)_* 
  \varphi_{ij} = (-1)^{n-j} f_1\wedge \cdots \wedge f_n \otimes e_1^* \wedge
\cdots e_{j-1}^* \wedge e_{j+1}^* \wedge \cdots \wedge e_n^* \wedge e_i^*. 
$$
Then the image of this tangent basis has dimension $m-n+1$. This happens at any point of 
$R(V,W)$. Now, the image of an application of constant rank is locally a submanifold.

Finally we have to check that the counterimages of $\bigwedge^n$ are connected, i.e.\
given two morphisms $\varphi_0$ and $\varphi_1$ such that $\bigwedge^n \varphi_0= 
\bigwedge^n \varphi_1$ then there exists a path $\{ \varphi_t\}_{t\in[0,1]}$ 
connecting
the two morphisms and satisfying $\bigwedge^n \varphi_t= \bigwedge^n \varphi_0$.
For this, note that the kernels of $\varphi_0$ and $\varphi_1$ coincide. Therefore there
exists an endomorphism $A$ in $GL(W)$ such that $A\varphi_0=\varphi_1$. Such $A$ is
forced to be in $SL(W)$. Now fix a path $A_t$, $t\in [0,1]$, connecting the identity with
$A$ and put $\varphi_t=A_t\varphi_0$.
\hfill $\Box$

This Lemma extends trivially to vector bundles to obtain the following
\begin{lemma} \label{beha_r}
Let $E,F$ be vector bundles of ranks $m$ and $n$ ($m\geq n$) respectively, then the fibration 
$R(E,F)$, given at any point $x\in M$ by $\bigwedge^n (\Hom(E_x,F_x))-\{ 0 \}$, has 
smooth fibers which are open complex submanifolds of 
$\Hom(\bigwedge^n E_x, \bigwedge^n F_x)$ of dimension $m-n+1$. Moreover 
$R(E,F)$ is invariant under multiplication by a never null complex-valued function, and so
given any point $d\in R(E,F)$ we have, using the standard identification between a
vector space and its tangent space at a point, that $d\in T_d R(E,F)$.
\end{lemma}

\subsubsection{Generalized asymptotically holomorphic sequences of sections of vector
bundles}
Now we recall the process of construction of a sequence of symplectic determinantal 
submanifolds. Let $E$, $F$ be vector bundles of ranks $r_e$ and $r_f$, respectively, and
suppose $r_e\geq r_f$. Write
$$
  r=\min\{r_e,r_f\}= r_f.
$$
Fix a generic section $s$ of the universal bundle $\SU_{ef}$ over $\Bi(r_e,r_f,N)$.
We embed $M$ in $\Bi(r_e,r_f,N)$ constructing an asymptotically holomorphic sequence
$\phi_k$ of embeddings. Using Lemma \ref{strata_bor} and Proposition \ref{strata_int}
we assure that the sequence is transverse to the holomorphic determinantal varieties 
defined by $s$ in $\Bi(r_e, r_f, N)$. We can define a sequence of sections of
the bundles $E^*\otimes F \otimes L^{\ox 2k}$ as
$$
 s_k= \phi_k^* s. 
$$ 
We consider now the connection $\hat{\nabla}_k$ defined on 
$E^*\otimes F\otimes L^{\otimes 2k}$
as the pull-back of the canonical one defined in $\SU_{ef}$. Also we consider in $M$
the sequence of metrics $\hat{g}_k$ defined as the pull-back through $\phi_k$ of the
standard metric on the bigrassmanian $\Bi(r_e, r_f, N)$. 
Then using properties 1 and 2 of Definition \ref{def:ah}, we obtain that the sequence
$s_k$ is asymptotically holomorphic with respect to the fixed complex structure $J$ in $M$,
computing the derivatives respect to $\hat{\nabla}_k$ and the norms respect to 
$\hat{g}_k$. Analogously taking the pull-back of the connection associated to 
$\bigwedge^{r} \pi_e^*(\SU_e) \otimes \bigwedge^{r} \pi_f^*(\SU_f)$, 
we obtain connections for the bundles $\bigwedge^{r} (E^*\ox L^{\ox k}) 
\otimes \bigwedge^{r} (F \otimes L^{\otimes k})$. 
Then the sequence $\bigwedge^r s_k$ is asymptotically $J$-holomorphic with respect to 
these connections and to the metric $\hat{g}_k$. 

Now we look for a condition to express when the sections $\bigwedge^r s_k$ are transversal 
in a certain sense. The key property is

\begin{lemma}
Let $E$ and $F$ be vector bundles with connections $\nabla^e$ and $\nabla^f$ respectively.
Suppose $s$ is a section of the bundle of morphisms $E^*\otimes F$ equipped with the connection
$\nabla^{ef}$ induced by $\nabla^e$ and $\nabla^f$. If $\bigwedge^r s(x)\neq 0$ 
at a point $x\in M$, then $\nabla^{ef} \bigwedge^r s(x)\in T_{\wedge^r s(x)} R(E_x,F_x)$.
\end{lemma}

{\bf Proof.}
To check this we only have to show that the following diagram is commutative
$$
\begin{array}{ccc}
\Omega^0(E^*\otimes F) & \stackrel{id \oplus \nabla^{ef}}{\to} & \Omega^0(E^*\otimes F)
\bigoplus \Omega^1(E^*\otimes F) \\
\downarrow \bigwedge^r & & \downarrow T\bigwedge^r \\
\Omega^0(\bigwedge^r(E^*)\otimes \bigwedge^r F) & \stackrel{id \oplus \nabla^{ef}}{\to} & 
\Omega^0(\bigwedge^r(E^*)\otimes \bigwedge^r F)
\bigoplus \Omega^1(\bigwedge^r(E^*)\otimes \bigwedge^r F).
\end{array}
$$
The map $T\bigwedge^r$ is defined as 
\begin{eqnarray*}
 T\bigwedge\nolimits^r : \Omega^0(E^*\otimes F)
  \oplus \Omega^1(E^*\otimes F) & \to & \Omega^0(\bigwedge\nolimits^r(E^*) \otimes 
  \bigwedge\nolimits^r F)
\oplus \Omega^1(\bigwedge\nolimits^r(E^*) \otimes \bigwedge\nolimits^r F) \\
(s_0,s_1) & \mapsto &  \left(s_0, \lim_{t\to 0} \frac{\bigwedge\nolimits^r(s_0+ts_1)}{t}\right). 
\end{eqnarray*}
To check this one fixes local frames in $E$ and $F$ and carries out
the computation explicitly.
\hfill $\Box$

Given a generic section $s$ of the bundle of morphisms $E^* \otimes F$ then we denote by
$D_{r-2}^{\epsilon}$ the $\epsilon$-neighborhood of the determinantal set $D_{r-2}$
associated to $s$.

\begin{definition}
Let $E$ and $F$ be vector bundles over $M$ of ranks $r_e$ and $r_f$ ($r_e\geq r_f$)
respectively. Put $r=r_f$.
We say that the section $s$ is $\eta$-$\wedge^r$-transverse to $0$, for some $\eta>0$, 
if for any $x\in M-D_{r-2}^{\eta}$ such
that $|\bigwedge^r s(x)|<\eta$ then the covariant derivative 
$\hat{s}(x)=\nabla \bigwedge^r s(x)$
has rank $r_e-r_f+1$ and also there exists a right inverse $\theta: T_{\wedge^r s(x)} R^r(E,F)
\to T_xM$ of $\hat{s}(x)$ with norm less that $\eta^{-1}$.
\end{definition}

We cannot impose the estimated transversality near the stratum $D_{r-2}$ because the
section $\bigwedge^r s$ is always critical in that stratum, so if we want to obtain a
notion of estimated transversality we need to remove a neighborhood 
of $D_{r-2}$.

Observe that given any small $\eta>0$,
the section $s$ is $\eta$-$\wedge^r$-transverse to $0$.

Using that $\phi_k(M)$ is transverse to $D_{r-1}$ we can check that $s_k$ is 
$\eta'$-$\wedge^r$-transverse
to $0$ on $M$, for some universal $\eta'>0$, with the connections and metrics defined in 
the prededent lines. Observe that to guarantee this property is absolutely necessary that 
the minimum distance from $\phi_k(M)$ to $D_{r-2}$ be greater than $\eta$, but this is true 
by construction.

\subsubsection{Proof of Proposition \ref{lefs}.}
We have as starting data a sequence of asymptotically holomorphic
sections of the bundles $E^*\otimes F\otimes L^{\ox 2k}$ obtained by pull-back
of a fixed section $s$ of the universal bundle $\SU_{ef}$. 
As before, we may suppose that $r_e\geq r_f$ and write $r=r_f$. Therefore the only
non-empty stratum is $D_{r-1}^k$, by assumption. We assume also
that $s_k$ is $\eta$-$\wedge^r$-transverse to $0$, for a universal $\eta>0$.
The stratum $D_{r-2}$ is empty and so the $\eta$-$\wedge^r$-transversality is
checked all over $M$. We can follow the ideas of 
\cite{Do96,Au97} to develop the proof. 

We define the function $f_k=\log |\bigwedge^r s_k|^2$. Clearly $f_k(-\infty)=
D_{r-1}^k$. Denote the complex dimension of $D_{r-1}^k$ by $N$. We are going to show that
all the critical points of $f_k$ are of index at least $N+1$. Therefore a standard
Morse-theoretic argument will finish the proof. 

Denote $\s_k=\bigwedge^r s_k$. First notice that if $x$ is a critical point of $|\s_k|^2$
then $\s_k(x)$ is not in the image of $\nabla \s_k$ and so $\nabla
\s_k$ is not surjective to $T_{\s_k(x)} R(E_x,F_x)$. It follows
from the $\eta$-$\wedge^r$-transversality property that $|\s_k(x)|>\eta$. 

Now we differentiate $f_k$ to obtain
$$ 
\partial f_k = \frac{1}{|\s_k|^2} ( \langle \partial \s_k, \s_k \rangle +
\langle \s_k, \bar{\partial} \s_k \rangle). 
$$
At a critical point $x$, $\partial f_k(x)=0$. Using the asymptotic holomorphic
bounds we obtain
\begin{equation}
|\langle \partial \s_k, \s_k \rangle|= |\langle \bar{\partial} \s_k,
\s_k \rangle|\leq  Ck^{-1/2}|\s_k|. \label{quiet}
\end{equation}
Differentiating a second time we obtain, evaluating at a critical point, the expression
$$ 
 \bar{\partial}\partial \log |\s|^2=\frac{1}{|\s|^2}(\langle \bar{\partial}\partial \s,
\s \rangle- \langle \partial \s, \partial \s \rangle
+ \langle \bar{\partial} \s, \bar{\partial} \s \rangle +
\langle \s, \partial \bar{\partial} \s \rangle), 
$$
where we omit the subindex $k$ for simplicity. Recall that $\bar{\partial}\partial+
\partial \bar{\partial}$ equals the $(1,1)$-part of the curvature of the bundle
$\bigwedge^r (E^*\otimes L^{\otimes k}) \otimes \bigwedge^r (F\otimes L^{\otimes k})$. 
Its $(1,1)$-curvature $R$ is the pull-back through $\phi_k$ of the 
$(1,1)$-curvature $\tilde R$ of $\bigwedge^r\SU_e \ox \bigwedge^r\SU_f$. 
So we obtain
$$ 
\bar{\partial} \partial f_k = \frac{1}{|\s|^2}(\langle R\s, \s \rangle -\langle
 \partial \bar{\partial} \s, \s \rangle + \langle \s, \partial \bar{\partial} 
 \s \rangle - \langle \partial \s, \partial \s \rangle +  
 \langle \bar{\partial} \s, \bar{\partial} \s \rangle). 
$$
We define the subspace
$$ 
 \CrV = \{ v\in T_xM \big| \nabla_v \s(x)= \lambda \s(x), 
  \text{ for some }\l \in \C\}. 
$$
Using the inequality \eqref{quiet} we obtain, for any $v\in \CrV$, that
$$ 
 |\la \partial_v \s, \s \ra |= |\partial_v \s| |\s|\leq C k^{-1/2}|\s|. 
$$
Restricting $\bbd\bd f_k$ to $\CrV$, it equals to 
$\frac{1}{|\s|^2}\la R\s, \s \ra
+O(k^{-1/2})$. Denote the Hessian of $f$ by $H_f$. 
We know that $H_f(u)+H_f(Ju)=-2i 
\bar{\partial} \partial f_k (u,Ju)=
-2i\frac{1}{|\s|^2}\la R(u,Ju) \s, \s \ra
+O(k^{-1/2})$, for any unit vector $u\in\CrV$. 
We claim that it is possible to bound above the expression 
\begin{equation}
   -2i\frac{1}{|\s|^2}\la R(u,Ju)\s, \s \ra \label{target}
\end{equation}
by a universal strictly negative constant, where $u$
is a unitary vector. For this we need to estimate the curvature $R$.
We start by computing the curvature of the universal bundle $\SU$ over 
the grassmannian $\Gr(r,N)$. We use the local expression of the 
curvature of $\SU^*$ from \cite[page 82]{We73},
$$ 
  R_{\SU^*}= h^{-1}\overline{df}^t \wedge df 
  -h^{-1}\overline{df}^t f h^{-1} \wedge \bar{f}^t df, 
$$
where $f=(f_1, \ldots, f_r)$ is a frame in an open neighborhood of 
$\Gr(r,N)$ and $h=\bar{f}^tf$. We may assume that we are at the point 
$\Pi_0=[\id|{\bf 0}]$ of the grassmannian, after suitable change of 
coordinates. Select the following holomorphic local frame,
$$
  f= ((1, 0,  \stackrel{(r-1)}{\ldots}, 0, z_{11}, \ldots, z_{1,n-r}), 
   \ldots, (0, \ldots, 0,1, z_{r1}, \ldots, z_{r,n-r})),
$$
So at the point $\Pi_0$ we obtain $R_{\SU^*}=\overline{df}^t\wedge df$ and
$$ 
 R_\SU = df^t\wedge \overline{df}.
$$
In the trivialization $(z_{jk})$ we take the standard basis 
$e_{jk}=\frac{\partial}{\partial z_{jk}}$. We obtain
$R_\SU (e_{jk}, ie_{jk})=-i b_{jj}$, where the endomorphism $b_{jj}$ is 
defined as $e_j\otimes e_j^*$. So the endomorphism $-iR_\SU (u,Ju)$ is 
semi-definite negative for $u \in T_{\Pi_0}\Gr(r,N)$ non-zero. This 
implies also that $-i R_{\bigwedge^k \SU}(u,Ju)$ is semi-definite 
negative, for $1\leq k\leq r$. 
Moreover computing $-i R_{\bigwedge^{r} \SU}(u,Ju)$ in 
$\Gr(r,N)$, or recalling that $\SU$ is very ample, we obtain that 
this endomorphism is definite negative. 
Returning to $\Bi(r_e,r_f,N)$ with $r=r_f\leq r_e$, we have that 
the curvature of $\bigwedge^r\SU_e \ox \bigwedge^r\SU_f$ is 
$$ 
  \tilde{R}= R_{\pi_e^* \bigwedge\nolimits^r \SU_e}\otimes \id_1 + \id_{\nu}
  \otimes R_{\pi_f^*\bigwedge\nolimits^r \SU_f}, 
$$
where $\nu=\left(\begin{array}{c} r_e\\ r \end{array}\right)$. 
So $\tilde{R}(u,Ju)$ is definite negative, for $u\in T\Bi(r_e,r_f,N)$ 
unitary vector.
Using that the sequence of embeddings $\phi_k=(\phi_k^e, \phi_k^f)$ satisfies 
properties $1$ and $2$ of Definition \ref{def:ah}, we get that the
expression \eqref{target} is bounded above by a universal strictly 
negative number.

Therefore, for any unitary $u\in \CrV$, 
$H_f(u)+H_f(Ju)$ is negative for $k$ large enough.
Recall that from the definition we obtain that $\dim \CrV\geq 2N+2$. 
Suppose that there exists a subspace $P\in T_xM$ of real dimension at 
least $2n-N$ such that $H_f$ in 
non-negative. The dimension of $P\bigcap JP$ is at least $2n-2N$, and there the function
$H_f(\cdot)+H_f(J\cdot)$ is, obviously, non-negative.   
Therefore $P\bigcap JP$ has to intersect trivially with $\CrV$ but $\dim P\bigcap JP+
\dim \CrV \geq 2n+2$, and this is clearly impossible. So such space $P$ does not exist
and then the index of $f_k$ at $x$ is greater than $N$. This finishes the proof.
\hfill $\Box$

\subsection{Chern classes of the constructed submanifolds} \label{chern}
For computing the Chern classes of determinantal submanifolds,
we shall use the results of Harris and Tu in \cite{HT84}.
All their results are stated for holomorphic determinantal
submanifolds in a holomorphic manifold, but they apply without 
the condition of integrability of the complex structure. We state
the formulas that we shall use. Following Subsection \ref{bigr}
we denote $r_e=\rank\,E$, $r_f=\rank\, F$, $2n=\dim M$ and
$D_r$ is the $r$-determinantal loci of a bundle map
$\varphi:E\to F$ constructed in Theorem \ref{main_deter}. 
First of all, set 
\begin{equation}
 \Delta_{i_1, \ldots, i_{r_e-r}} = \left| \begin{array}{cccc}
 c_{r_f-r+i_1} & c_{r_f-r+i_1+1} & \cdots & \\
 c_{r_f-r+i_2-1} & c_{r_f-r+i_2} & \cdots & \\
 & & \ddots & \\
 & & \cdots & c_{r_f-r+i_{r_e-r}}
 \end{array}
 \right| ,
 \label{delta}
\end{equation}
where $c_j=c_j(F-E)$. For instance, $\Delta_{0, \ldots, 0}=\Delta
=\PD([D_r])$, which is the classical Porteous formula for the 
homology class of a determinantal locus. We can suppose that the 
indices $i_j$ are decreasing, and so if we have any index $i_j=0$
we do not write it, e.g.\ $\Delta_{2,1,0}=\Delta_{2,1}$.

In \cite{HT84} a complete description of the Chern numbers of the tangent
bundle of a determinantal submanifold is performed, supposing that
$D_{r-1}= \emptyset$ and so $D_r$ is smooth. We concentrate ourselves
in the cases $\dimC D_r=1$ and $\dimC D_r=2$, where Harris and Tu obtain
the following formulas:
\begin{enumerate}
\item For $\dimC M=(r_e-r)(r_f-r)+1$, then $\dimC D_r=1$. We have
  $$
    n_1(D_r)= \la c_1(D_r), [D_r] \ra =(c_1(M)+(r_e-r)c_1(E- F))\Delta + 
    (r_e-r_f)\Delta_1.
  $$
\item For $\dimC M=(r_e-r)(r_f-r)+2$, then $\dimC D_r=2$. We have
\begin{eqnarray*}
 n_{11}(D_r) & = & \langle c_1^2(D_r), [D_r] \rangle =
    (c_1(M)+(r_e-r)c_1(E- F))^2\cdot \Delta+ \\
  &+ & 2(r_e-r_f)(c_1(M)+(r_e-r)c_1(E- F))\cdot \Delta_1 +(r_e-r_f)^2(\Delta_2
    +\Delta_{11}), \\
 n_2(D_r)&=&\la c_2(D_r), [D_r] \ra =\big( c_2(M) +(r_e-r)c_1(M)c_1(E- F)+ \\
  &+ & (r_e-r)(c_2(E)-c_2(F))+\binom{r_e-r}{2}c_1^2(E)-(r_e-r)^2c_1(E)c_1(F)+\\
  &+ & \binom{r_e-r+1}{2} c_1^2(F) \big) \Delta+ \\
  &+ & ((r_e-r)c_1(M)+((r_e-r)(r_e-r_f)-1)c_1(E-F)) \Delta_1 \\
  &+ & \frac{1}{2} ((r_e-r_f)^2+(r_e-r)+(r_f-r)-2)\Delta_2+ \\
  &+ & \frac{1}{2} ((r_e-r_f)^2-(r_e-r)-(r_f-r)-2) \Delta_{11}. 
\end{eqnarray*}
\end{enumerate}

In our case, we are going to apply the above formulas to morphisms
$\varphi:E\ox (L^*)^{\ox k} \to F \ox L^{\ox k}$. We have the following
asymptotic expansions for Chern classes (we write $\o_k=\frac{k\o}{2\pi}$
for simplicity)
\begin{eqnarray} 
 & & c_p(F\ox L^{\ox k}) = \binom{r_f}{p} \o_k^p +O(k^{p-1}), \nonumber\\
 & & c_p(E\ox (L^*)^{\ox k}) =\binom{r_e}{p} (-\o_k)^p + O(k^{p-1}),
  \nonumber \\
 & & c_p=c_p(F\ox L^{\ox k}-E\ox (L^*)^{\ox k})= 
  \text{Coeff}_{x^p} \frac{(1+x)^{r_f}}{(1-x)^{r_e}} \o_k^p +O(k^{p-1})=
   \nonumber \\
 & & \qquad = \sum_{i=0}^{r_f} \binom{r_f}{i} \binom{r_e+p-i-1}{p-i} 
    \o_k^p +O(k^{p-1}). \label{cp}
\end{eqnarray}

We are going to give two families of examples to show that the symplectic manifolds 
obtained here are more general than those in \cite{Au97}.

\subsubsection{Example 1.}
Choose $\dimC M=(r_e-r)(r_f-r)+1$ and so we can apply the formulas for the
complex $1$-dimensional case. Also suppose that $r=1$ and $r_e=2$, so 
$\dimC M= r_f=n>1$. By Proposition \ref{lefs} the submanifolds $D_1$
are connected. Now $\PD[D_1]=\D=c_{n-1}$ and $\D_1=c_n$. Using
\eqref{cp} we get that
\begin{eqnarray*}
  \vol_{\o_k}(D_1) &=&\D\o_k=(n2^{n-1} +O(k^{-1}))\vol_{\omega_k} (M), \\
  n_1(D_1) &=& -(n+2)\o_k\D + (2-n)\D_1 +O(k^{-1})\vol_{\omega_k} (M) =\\
       &=&  (-(n+2)n2^{n-1} + (2-n)( n 2^{n-1}+2^n)+ O(k^{-1}))\vol_{\omega_k} (M) \\
  \frac{n_1(D_1)}{\vol_{\o_k}(D_1)}&=& -2-2n+\frac{4}{n}+O(k^{-1}).
\end{eqnarray*}
 To compare with the Auroux' case we compute the precedent symplectic 
 invariants for this situation. Denote by $Z$ the zero set of a transverse 
 section of a bundle of the form $E\otimes L^{\ox k}$, we choose 
 $\rank\, E= n-1$ to set up the comparison. Suppose that $Z$ is symplectic. 
 Using Proposition 5 in \cite{Au97} we obtain
\begin{eqnarray*}
 \vol_{\o_k}(Z) & = & (1+O(k^{-1})) \vol_{\omega_k} (M), \\
 n_1(Z) & = & (1-n+O(k^{-1})) \vol_{\omega_k} (M), \\
 \frac{n_1(Z)}{\vol_{\omega_k}(Z)} & = & 1-n+O(k^{-1}).
\end{eqnarray*}
 Therefore there does not exist any $n\geq 2$ such that the
 quotients $\frac{n_1(D_1)}{\vol_{\omega_k}(D_1)}$ coincide with the quotients
 $\frac{n_1(Z)}{\vol_{\omega_k}(Z)}$,
 obviously for $k$ large enough. So Auroux' sequences of submanifold 
 are not symplectomorphic to our sequences of determinantal submanifolds.

 To check that, for $k$ large, our determinantal submanifolds do not coincide with 
 Auroux' examples we work as follows. Suppose that for integers $k_1, k_2$ the
 submanifold $D_1=D_1^{k_1}$ is isotopic to $Z=Z_{k_2}$. Then they define the same
 cohomology class and hence $n2^{n-1}k_1=k_2+O(1)$. Also $n_1(D_1)=n_1(Z)$ implies
 $(-2-2n+\frac 4n)k_1=(1-n)k_2 +O(1)$. So, for large enough $k$'s, $(1-n)n2^{n-1}=
 -2-2n+\frac 4n$ and hence $n=2$. Therefore for $n>2$ and large $k$ we get new
 examples of symplectic submanifolds.

Note that for $n=r_e=r_f=2$, the determinantal set $D_1$ for a morphism
$\varphi:E\ox (L^*)^{\ox k} \to F \ox L^{\ox k}$ is 
the zero set of the section $\bigwedge^2 \varphi$ of 
$\bigwedge\nolimits^2 E^*\ox \bigwedge\nolimits^2 F \ox L^{\ox 4k}$. Since
this zero set is smooth of the expected codimension, our example is just
one of Auroux' examples.

\subsubsection{Example 2.}
Now, choose $\dimC M=(r_e-r)(r_f-r)+2$ and so we can apply the formulas for 
the complex $2$-dimensional case. Again we suppose that $r=1$ and $r_e=2$, so
$\dimC M= r_f+1=n>2$.  By Proposition \ref{lefs} these submanifolds
are connected. In this case we have
\begin{eqnarray*}
 \vol_{\omega_k}(D_1) & = & ((n-1)2^{n-2}+O(k^{-1})) \vol_{\o_k} (M), \\
 n_{11}(D_1) & = & (4(n-1)(n^2-5)2^{n-2} +O(k^{-1}))  \vol_{\o_k} (M)\\
 n_2(D_1) & = & (2(n^2+n-4)(n-1)2^{n-2} +O(k^{-1}) \vol_{\o_k} (M) \\
 \frac{n_2(D_1)}{n_{11}(D_1)} &=& \frac{n^2+n-4}{2(n^2-5)} +O(k^{-1}).
\end{eqnarray*}
For the Auroux' case with $\rank\, E=n-2$ we obtain
\begin{eqnarray*}
 \vol_{\omega_k}(Z) & = & (1+ O(k^{-1})) \vol_{\omega_k} (M), \\
 n_{11}(Z) & = & ((n-2)^2+ O(k^{-1})) \vol_{\omega_k} (M), \\
 n_2(Z) & = & \left(\frac{(n-1)(n-2)}{2}+ O(k^{-1})\right) \vol_{\o_k} (M) \\
 \frac{n_2(Z)}{n_{11}(Z)} &=& \frac{n-1}{2(n-2)} +O(k^{-1}).
\end{eqnarray*}
If we compute the symplectic invariants $\frac{n_{11}(Z)}{\vol_{\omega_k}(Z)}$
and $\frac{n_2(Z)}{\vol_{\omega_k}(Z)}$, it is easy to verify that Auroux'
submanifolds are not symplectomorphic to the determinantal ones constructed in
this example.

Moreover, for $4$-manifolds, the numbers $n_2=\chi$ and $n_{11}=(2\chi+3\sigma)/4$ are
topological invariants.
Therefore $\frac{n_2}{n_{11}}$ is a topological invariant. Comparing
the Auroux' case and the determinantal example we find that these symplectic 
submanifolds are not even {\it homeomorphic}, for $k$ large enough (even choosing
different $k$'s in either case).

In general, it is clear that the determinantal class is quite bigger than the Donaldson-Auroux
one. We could compute more examples and more precise invariants using recent results
from algebraic geometry about the topology of determinantal submanifolds. 
As a reference it could be useful \cite{HT84b, Pr88, PP91}. Remark that in these
references the computations are performed even in the singular case. To adapt them
to the symplectic category we would need to define the Segre classes of a singular
symplectic manifold. This definition seems quite natural.

\end{document}